\RequirePackage{fix-cm}
\documentclass{svjour3}                      
\smartqed  
\usepackage{times}
\usepackage{latexsym}
\usepackage{cite}
\usepackage[misc,geometry]{ifsym}
\usepackage{amsmath} 	
\usepackage{amssymb}  	
\usepackage{amsfonts}
\usepackage{graphicx}
\usepackage{subfigure}
\usepackage{epstopdf}
\usepackage{dblfloatfix}
\usepackage{tikz}
\usetikzlibrary{shapes,matrix,decorations.shapes}

\usepackage{array,xcolor}
\usepackage{algorithm}
\usepackage{algpseudocode}
\usepackage{multirow}
\usepackage{threeparttable}
\usepackage{booktabs}
\usepackage[misc]{ifsym}

\usepackage{enumerate}

\usepackage{textgreek}

\usepackage[colorlinks=true,      
			linkcolor=black,      
			citecolor=black,      
			filecolor=black,      
			urlcolor=blue ]{hyperref}

\DeclareMathOperator{\vect}{vec}
\DeclareMathOperator{\mat}{mat}

\newcommand{\algsize}{\small}
\newcommand{\probname}[1]{{\tt\small #1}}
\newcommand{\CDCS}{\mbox{CDCS}}

\begin{document}

\title{
Chordal decomposition in operator-splitting methods for sparse semidefinite programs
\thanks{YZ
and
GF
contributed equally. A preliminary version of part of this work appeared in~\cite{ZFPGWpd2016,ZFPGWhsde2016}.
 YZ is supported by Clarendon Scholarship and Jason Hu Scholarship.
GF was supported by EPSRC grant EP/J010537/1 {and by and EPSRC Doctoral Prize Fellowship}.
AP was supported in part by EPSRC Grant EP/J010537/1 and EP/M002454/1.
}}
\titlerunning{Chordal decomposition in operator-splitting methods for sparse SDPs}        

\author{Yang Zheng$^{1}$\and Giovanni Fantuzzi$^{2}$ \and \\ Antonis Papachristodoulou$^{1}$ \and Paul Goulart$^{1}$ \and Andrew~Wynn$^{2}$ }

\authorrunning{Y. Zheng, G. Fantuzzi, A. Papachristodoulou, P. Goulart and A. Wynn} 

\institute{
Y. Zheng (\Letter)\\
Tel.: +44-07511784230\\
\email{yang.zheng@eng.ox.ac.uk} \vspace{2 pt}\\
G. Fantuzzi \\
\email{giovanni.fantuzzi10@imperial.ac.uk}\vspace{2 pt}\\
A. Papachristodoulou \\
\email{antonis@eng.ox.ac.uk} \vspace{2 pt}\\
P. Goulart \\
\email{paul.goulart@eng.ox.ac.uk} \vspace{2 pt}\\
A. Wynn \\
\email{a.wynn@imperial.ac.uk} \vspace{2 pt}\\
$^{1}\;$Department of Engineering Science, University of Oxford, Parks Road, Oxford, OX1 3PJ, U.K. \\ 
$^{2} \;$Department of Aeronautics, Imperial College London, South Kensington Campus, SW7 2AZ, U.K.
}

\date{Received: date / Accepted: date}

\maketitle
\begin{abstract}

We employ chordal decomposition to reformulate a large and sparse semidefinite program (SDP), either in primal or dual standard form, into an equivalent SDP with smaller positive semidefinite (PSD) constraints. In contrast to previous approaches, the decomposed SDP is suitable for the application of first-order operator-splitting methods, enabling the development of efficient and scalable algorithms. In particular, we apply the alternating direction method of multipliers (ADMM) to solve decomposed primal- and dual-standard-form SDPs. Each iteration of such ADMM algorithms requires a projection onto an affine subspace, and a set of projections onto small PSD cones that can be computed in parallel. We also formulate the homogeneous self-dual embedding (HSDE) of a primal-dual pair of decomposed SDPs, and extend a recent ADMM-based algorithm to exploit the structure of our HSDE. The resulting HSDE algorithm has the same leading-order computational cost as those for the primal or dual problems only, with the advantage of being able to identify infeasible problems and produce an infeasibility certificate. All algorithms are implemented in the open-source MATLAB solver {{\CDCS}}. Numerical experiments on a range of large-scale SDPs demonstrate the computational advantages of the proposed methods compared to common state-of-the-art solvers.

\end{abstract}

\section{Introduction} \label{intro}

Semidefinite programs (SDPs) are convex optimization problems over the cone of positive semidefinite (PSD) matrices. Given $ b\in \mathbb{R}^m$, $C\in \mathbb{S}^n$, and matrices $A_1,\,\ldots,\,A_m \in \mathbb{S}^n$, the standard \textit{primal form} of an SDP is
\begin{equation}
\label{E:PrimalSDP}
    \begin{aligned}
    \min_{X} \quad & \langle C, X \rangle \\
    \text{subject to} \quad & \langle A_i, X\rangle = b_i, \quad i = 1, \ldots, m, \\
        & X \in \mathbb{S}^n_{+},
    \end{aligned}
\end{equation}
while the standard \textit{dual form} is
\begin{equation}
\label{E:DualSDP}
    \begin{aligned}
        \max_{y, \;Z} \quad & \langle b,y\rangle \\
        \text{subject to} \quad  & Z + \sum_{i=1}^m A_i\,y_i = C,\\
        & Z \in \mathbb{S}^n_{+}.
    \end{aligned}\\
\end{equation}
In the above and throughout this work, $\mathbb{R}^m$ is the usual $m$-dimensional Euclidean space, $\mathbb{S}^n$ is the space of $n \times n$ symmetric matrices, $\mathbb{S}^n_{+}$ is the cone of PSD matrices, and $\langle \cdot,\cdot \rangle$ denotes the inner product in the appropriate space, \emph{i.e.}, $\langle x,y\rangle = x^Ty$ for $x,y\in \mathbb{R}^m$ and $\langle X, Y\rangle = \mathrm{trace}(XY)$ for $X,Y\in \mathbb{S}^n$.
SDPs have found applications in a wide range of fields, such as control theory, machine learning, combinatorics, and operations research~\cite{boyd1994linear}. Semidefinite programming encompasses other common types of optimization problems, including linear, quadratic, and second-order cone programs~\cite{boyd2004convex}. Furthermore, many nonlinear convex constraints admit SDP relaxations that work well in practice~\cite{vandenberghe1996semidefinite}.

It is well-known that small and medium-sized SDPs can be solved up to any arbitrary precision in polynomial time~\cite{vandenberghe1996semidefinite} using efficient second-order interior-point methods (IPMs)~\cite{alizadeh1998primal, helmberg1996interior}. However, many problems of practical interest are too large to be addressed by the current state-of-the-art interior-point algorithms, largely due to the need to compute, store, and factorize an $m\times m$ matrix at each iteration.

A common strategy to address this shortcoming is to abandon IPMs in favour of simpler first-order methods (FOMs), at the expense of reducing the accuracy of the solution. For instance, Malick \emph{et al.} introduced regularization methods to solve SDPs based on a dual augmented Lagrangian~\cite{malick2009regularization}. Wen \emph{et al.} proposed an alternating direction augmented Lagrangian method for large-scale SDPs in the dual standard form~\cite{wen2010alternating}. Zhao \emph{et al.} presented an augmented Lagrangian dual approach combined with the conjugate gradient method to solve large-scale SDPs~\cite{zhao2010newton}. More recently, O'Donoghue \emph{et al.} developed a first-order operator-splitting method to solve the homogeneous self-dual embedding (HSDE) of a primal-dual pair of conic programs~\cite{ODonoghue2016}. The algorithm, implemented in the C package SCS \cite{scs}, has the advantage of providing certificates of primal or dual infeasibility.

A second major approach to resolve the aforementioned scalability issues is based on the observation that the large-scale SDPs encountered in applications are often structured and/or sparse~\cite{boyd1994linear}. Exploiting sparsity in SDPs is an active and challenging area of research~\cite{andersen2011interior}, with one main difficulty being that the optimal (primal) solution is typically dense even when the problem data are sparse. Nonetheless, if the {\it aggregate sparsity pattern} of the data is \textit{chordal} (or has sparse \textit{chordal extensions}), one can replace the original, large PSD constraint with a set of PSD constraints on smaller matrices, coupled by additional equality constraints~\cite{agler1988positive,grone1984positive,griewank1984existence,kakimura2010direct}. Having reduced the size of the semidefinite variables, the converted SDP can in some cases be solved more efficiently than the original problem using standard IPMs. These ideas underly the \emph{domain-space} and the \emph{range-space} conversion techniques in~\cite{fukuda2001exploiting,kim2011exploiting}, implemented in the MATLAB package SparseCoLO~\cite{fujisawa2009user}.

The problem with such decomposition techniques, however, is that the addition of equality constraints to an SDP often offsets the benefit of working with smaller semidefinite cones. One possible solution is to exploit the properties of chordal sparsity patterns directly in the IPMs: Fukuda \emph{et al.} used a positive definite completion theorem~\cite{grone1984positive} to develop a primal-dual path-following method~\cite{fukuda2001exploiting}; Burer proposed a nonsymmetric primal-dual IPM using Cholesky factors of the dual variable $Z$ and maximum determinant completion of the primal variable $X$~\cite{burer2003semidefinite}; and Andersen \emph{et al.} developed fast recursive algorithms to evaluate the function values and derivatives of the barrier functions for SDPs with chordal sparsity~\cite{andersen2010implementation}. Another attractive option is to solve the sparse SDP using FOMs: Sun \emph{et al.} proposed a first-order splitting algorithm for partially decomposable conic programs, including SDPs with chordal sparsity~\cite{sun2014decomposition}; Kalbat \& Lavaei applied a first-order operator-splitting method to solve a special class of SDPs with fully decomposable constraints~\cite{Kalbat2015Fast}; Madani \emph{et al.} developed a highly-parallelizable first-order algorithm for sparse SDPs with inequality constraints, with applications to optimal power flow problems~\cite{Madani2015ADMM};  Dall'Anese \emph{et al.} exploited chordal sparsity to solve SDPs with separable constraints using a distributed FOM~\cite{dall2013distributed}; {finally, Sun and Vandenberghe introduced several proximal splitting and decomposition algorithms for sparse matrix nearness problems involving no explicit equality constraints~\cite{sun2015decomposition}.}

In this work, we embrace the spirit of~\cite{ODonoghue2016, sun2014decomposition, Kalbat2015Fast, Madani2015ADMM, dall2013distributed,sun2015decomposition} and exploit sparsity in SDPs using a first-order  operator-splitting method known as the \emph{alternating direction method of multipliers} (ADMM). Introduced in the mid-1970s~\cite{glowinski1975approximation,gabay1976dual}, ADMM is related to other FOMs such as dual decomposition and the method of multipliers, and it has recently found applications in many areas, including covariance selection, signal processing, resource allocation, and classification; see~\cite{boyd2011distributed} for a review. In contrast to the approach in~\cite{sun2014decomposition}, which requires the solution of a quadratic SDP at each iteration, our approach relies entirely on first-order methods. Moreover, our ADMM-based algorithm works for generic SDPs with chordal sparsity and has the ability {to detect} infeasibility, which are key advantages compared to the algorithms in~\cite{Kalbat2015Fast,Madani2015ADMM,dall2013distributed,sun2015decomposition}.
More precisely, our contributions are:

\begin{enumerate}
  \item
  We apply two chordal decomposition theorems~\cite{grone1984positive,agler1988positive} to formulate \emph{domain-space} and \emph{range-space} conversion frameworks for the application of FOMs to standard-form SDPs with chordal sparsity. These are analogous to the conversion methods developed in~\cite{fukuda2001exploiting, kim2011exploiting} for IPMs, but we introduce two sets of slack variables that allow for the separation of the conic and the affine constraints when using operator-splitting algorithms. To the best of our knowledge, this extension has never been presented before, and its significant potential is demonstrated in this work.

\vspace{0.25em}

\item We apply ADMM to solve the domain- and range-space converted SDPs, and show that the resulting iterates of the ADMM algorithms are the same up to scaling. The iterations are computationally inexpensive: the positive semidefinite (PSD) constraint is enforced via parallel projections onto small PSD cones---a much more economical strategy than that in~\cite{sun2014decomposition}---while imposing the affine constraints requires solving a linear system with constant coefficient matrix, the factorization/inverse of which can be cached before iterating the algorithm. Note that the idea of enforcing a large sparse PSD constraint by projection onto multiple smaller ones has also been exploited in~\cite{Madani2015ADMM, dall2013distributed} in the special context of optimal power flow problems { and in~\cite{sun2015decomposition} for matrix nearness problems}.

    \vspace{0.25em}

  \item  We formulate the HSDE of a converted primal-dual pair of sparse SDPs. In contrast to~\cite{sun2014decomposition, dall2013distributed, Kalbat2015Fast, Madani2015ADMM}, this allows us to compute either primal and dual optimal points, or a certificate of infeasibility. {We then extend the algorithm proposed in~\cite{ODonoghue2016}, showing that the structure of our HSDE can be exploited to solve a large linear system of equations extremely efficiently through a sequence of block eliminations.} As a result, we obtain an algorithm that is more efficient than the method of \cite{ODonoghue2016}, irrespectively of whether this is used on the original primal-dual pair of SDPs (before decomposition) or on the converted problems. In the former case, the advantage comes from the application of chordal decomposition to replace a large PSD cone with a set of smaller ones. In the latter case, efficiency is gained by the proposed sequence of block eliminations.

\vspace{0.25em}

\item We present the MATLAB solver {{\CDCS}} (Cone Decomposition Conic Solver), which implements our ADMM algorithms. {{\CDCS}} is the first open-source first-order solver that exploits chordal decomposition and can detect infeasible problems. We test our implementation on large-scale sparse problems in SDPLIB \cite{borchers1999sdplib}, selected sparse SDPs with nonchordal sparsity pattern~\cite{andersen2010implementation}, and randomly generated SDPs with block-arrow sparsity patterns~\cite{sun2014decomposition}. The results demonstrate the efficiency of our algorithms compared to the interior-point solvers SeDuMi \cite{sturm1999using} and the first-order solver SCS \cite{scs}.
\end{enumerate}

The rest of the paper is organized as follows. Section~\ref{sec:preliminaries} reviews chordal decomposition and the basic ADMM algorithm. Section~\ref{sec:decompositionSDPs} introduces our conversion framework for sparse SDPs based on chordal decomposition. We show how to apply the ADMM to exploit domain-space and range-space sparsity in primal and dual SDPs in Section~\ref{sec:primal&dual}. Section~\ref{sec:hsde} discusses the ADMM algorithm for the HSDE of SDPs with chordal sparsity. The computational complexity of our algorithms in terms of floating-point operations is discussed in Section~\ref{Section:Complexity}. {{\CDCS}} and our numerical experiments are presented in Section~\ref{sec:simulation}. Section~\ref{sec:conclusion} concludes the paper.

\section{Preliminaries}
\label{sec:preliminaries}

\subsection{A review of graph theoretic notions}

We start by briefly reviewing some key graph theoretic concepts (see~\cite{godsil2013algebraic, blair1993introduction} for more details). A graph $\mathcal{G}(\mathcal{V},\mathcal{E})$ is defined by a set of vertices $\mathcal{V}=\{1,2,\ldots,n\}$ and a set of edges $\mathcal{E} \subseteq \mathcal{V} \times \mathcal{V}$. A graph $\mathcal{G}(\mathcal{V},\mathcal{E})$ is called {\it complete} if any two nodes are connected by an edge. A subset of vertices $\mathcal{C}\subseteq \mathcal{V} $ such that $ (i,j) \in \mathcal{E}$ for any distinct vertices $ i,j \in \mathcal{C} $, \emph{i.e.}, such that the subgraph induced by $\mathcal{C}$ is complete,
is called a \emph{clique}. The number of vertices in $\mathcal{C}$ is denoted by $\vert \mathcal{C} \vert$. If $\mathcal{C}$ is not a subset of any other clique, then it is referred to as a \emph{maximal clique}. A \emph{cycle} of length $ k $ in a graph $\mathcal{G}$ is a set of pairwise distinct nodes $ \{v_1,v_2,\ldots,v_k\}\subset \mathcal{V} $ such that $ (v_k,v_1) \in \mathcal{E} $ and $ (v_i,v_{i+1}) \in \mathcal{E} $ for $ i=1,\ldots,k-1 $. A \emph{chord} is an edge joining two non-adjacent nodes in a cycle. A graph $\mathcal{G}$ is undirected if $(v_i,v_j) \in \mathcal{E} \Leftrightarrow (v_j,v_i) \in \mathcal{E}$.

{An undirected graph} $\mathcal{G}$ is called \emph{chordal} (or \emph{triangulated}, or a \emph{rigid circuit}~\cite{vandenberghe2014chordal}) if every cycle of length greater than or equal to four has at least one chord. Chordal graphs include several other classes of graphs, such as acyclic undirected graphs (including trees) and complete graphs. Algorithms such as the maximum cardinality search~\cite{tarjan1984simple} can test chordality and identify the maximal cliques of a chordal graph efficiently, {\it i.e.}, in linear time in terms of the number of nodes and edges.
Non-chordal graphs can always be \emph{chordal extended}, \emph{i.e.}, extended to a chordal graph, by adding additional edges to the original graph. Computing the chordal extension with the minimum number of additional edges is an NP-complete problem~\cite{yannakakis1981computing}, but several heuristics exist to find good chordal extensions efficiently~\cite{vandenberghe2014chordal}.

\begin{figure}[t]
    \centering
    \setlength{\abovecaptionskip}{0em}
    \setlength{\belowcaptionskip}{0em}
    \subfigure[]
    { \label{fig:e}
    	  \includegraphics[scale=1]{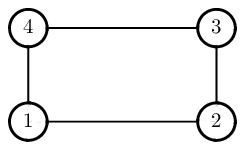}
    }
    \hspace{50pt}
    \subfigure[]
    { \label{fig:f}
      \includegraphics[scale=1]{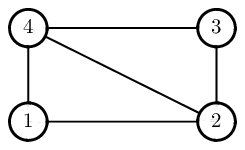}
    }
    \caption{(a) Nonchordal graph: the cycle (1-2-3-4) is of length four but has no chords. (b) Chordal graph: all cycles of length no less than four have a chord; the maximal cliques are $\mathcal{C}_1 = \{1,2,4\}$ and $\mathcal{C}_2 = \{2,3,4\}$.}
    \label{F:ChordalGraph}
\end{figure}

Fig.~\ref{F:ChordalGraph} illustrates these concepts. The graph in Fig.~\ref{F:ChordalGraph}(a) is not chordal, but can be chordal extended to the graph in Fig.~\ref{F:ChordalGraph}(b) by adding the edge $(2,4)$. The chordal graph in Fig.~\ref{F:ChordalGraph}(b) has two maximal cliques, $ \mathcal{C}_1 = \{1,2,4\}$ and $\mathcal{C}_2 = \{2,3,4\}$. Other examples of chordal graphs are given in Fig.~\ref{F:ExampleChordal}. 

\subsection{Sparse matrix cones and chordal decomposition}

\begin{figure}[t]
    \centering
    \setlength{\abovecaptionskip}{0em}
    \setlength{\belowcaptionskip}{0em}
    \newcommand{\fighspace}{\hspace{1cm}}
    \subfigure[]
    { \label{fig:d}
      \includegraphics[scale=0.9]{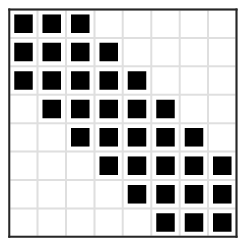}
    } \fighspace
    \subfigure[]
    { \label{fig:g}
      \includegraphics[scale=0.9]{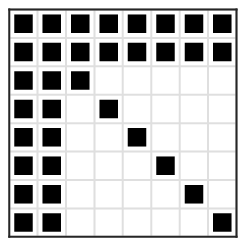}
    } \fighspace
    \subfigure[]
    { \label{fig:h}
      \includegraphics[scale=0.9]{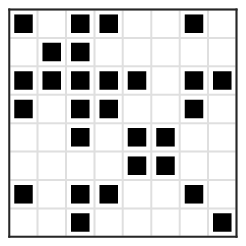}
    }
    \caption{{Sparsity patterns of $8\times 8$ matrices: (a) banded sparsity pattern; (b) ``block-arrow'' sparsity pattern; (c) a generic sparsity pattern.}}
    \label{F:ExampleSparsityPattern}
\end{figure}

\begin{figure}
    \centering
    \setlength{\abovecaptionskip}{0em}
    \setlength{\belowcaptionskip}{0em}
    \newcommand{\fighspace}{\hspace{1cm}}
    \subfigure[]
    { \label{fig:a}
      \includegraphics[scale=1]{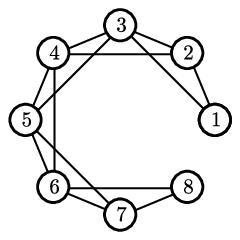}
    } \fighspace
    \subfigure[]
    { \label{fig:b}
      \includegraphics[scale=1]{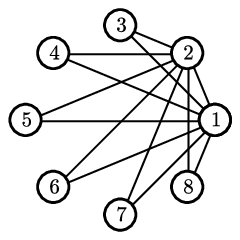}
    } \fighspace
    \subfigure[]
    { \label{fig:c}
      \includegraphics[scale=1]{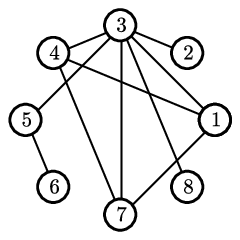}
    }
    \caption{{Graph representation of the matrix sparsity patterns illustrated in Fig.~\ref{F:ExampleSparsityPattern}(a)--(c), respectively.}}
    \label{F:ExampleChordal}
\end{figure}

The sparsity pattern of a symmetric matrix $X \in \mathbb{S}^n$ can be represented by an undirected graph $\mathcal{G}(\mathcal{V},\mathcal{E})$, and vice-versa. For example, the sparsity patterns illustrated in Fig.~\ref{F:ExampleSparsityPattern} correspond to the graphs in Fig.~\ref{F:ExampleChordal}. With a slight abuse of terminology, we refer to the graph $\mathcal{G}$ as the sparsity pattern of $X$. Given a clique $\mathcal{C}_k$ of $\mathcal{G}$, we define a matrix $E_{\mathcal{C}_k} \in \mathbb{R}^{\vert  \mathcal{C}_k\vert  \times n}$ as
$$
    (E_{\mathcal{C}_k})_{ij} = \begin{cases} 1, \quad \text{if } {\mathcal{C}_k}(i) = j \\ 0, \quad \text{otherwise} \end{cases}
$$
where $\mathcal{C}_k(i)$  is the $i$-th vertex in $\mathcal{C}_k$, sorted in the natural ordering. Given $X \in \mathbb{S}^n$, the matrix $E_{\mathcal{C}_k}$ can be used to select a principal sub-matrix defined by the clique $\mathcal{C}_k$, \emph{i.e.}, $  E_{\mathcal{C}_k}XE_{\mathcal{C}_k}^T \in \mathbb{S}^{\vert  \mathcal{C}_k\vert }$. In addition, the operation $E_{\mathcal{C}_k}^TYE_{\mathcal{C}_k}$ creates an $n \times n$ symmetric matrix from a $\vert\mathcal{C}_k\vert \times \vert \mathcal{C}_k\vert $ matrix. For example, the chordal graph in Fig.~\ref{F:ChordalGraph}(b) has a maximal clique $\mathcal{C}_1 = \{1,2,4\}$, and for {$X\in\mathbb{S}^4$} and $Y\in\mathbb{S}^3$ we have
$$
E_{\mathcal{C}_1} = \begin{bmatrix} 1 & 0 & 0& 0 \\ 0 & 1 & 0 & 0 \\ 0 & 0& 0 & 1\end{bmatrix},
\,\,
E_{\mathcal{C}_1}XE_{\mathcal{C}_1}^T = \begin{bmatrix} X_{11} & X_{12} & X_{14} \\ X_{21} & X_{22} & X_{24} \\
    X_{41} & X_{42} & X_{44}\end{bmatrix},
\,\,
E_{\mathcal{C}_1}^TYE_{\mathcal{C}_1}= \begin{bmatrix} Y_{11} & Y_{12} & 0 & Y_{13} \\ Y_{21} & Y_{22} & 0 & Y_{23}\\ 0 & 0 & 0 & 0 \\
    Y_{31} & Y_{32} & 0 & Y_{33}\end{bmatrix}.
$$

Given an undirected graph $\mathcal{G}(\mathcal{V}, \mathcal{E})$, let $\mathcal{E}^* = \mathcal{E} \cup \{(i,i),\, i \in \mathcal{V}\}$ be a set of edges that includes all self-loops.
We define the space of sparse symmetric matrices {with sparsity pattern} $\mathcal{G}$ as
\begin{equation*}
\mathbb{S}^n(\mathcal{E},0) := \{ X \in \mathbb{S}^n :\;  X_{ij} =X_{ji} = 0 \text{ if } (i,j) \notin \mathcal{E}^* \},
\end{equation*}
and the cone of sparse PSD matrices as
\begin{equation*}
\mathbb{S}^n_+(\mathcal{E},0) := \{ X \in \mathbb{S}^n(\mathcal{E},0):\, X\succeq 0  \},
\end{equation*}
{where the notation $X\succeq 0$ indicates that $X$ is PSD.}
Moreover, we consider the cone
\begin{equation*}
\mathbb{S}^n_+(\mathcal{E},?) := \mathbb{P}_{\mathbb{S}^n(\mathcal{E},0) }( \mathbb{S}^n_+ )
\end{equation*}
given by the projection of the PSD cone onto the space of sparse matrices $\mathbb{S}^n(\mathcal{E},0)$ with respect to the usual Frobenius matrix norm (this is the norm induced by the usual trace inner product on the space of symmetric matrices). It is not difficult to see that $X\in \mathbb{S}^n_+(\mathcal{E},?) $ if and only if it has a positive semidefinite completion, {\it i.e.}, if there exists {a PSD matrix $M$} such that $M_{ij}=X_{ij}$ when $(i,j)\in \mathcal{E}^*$.

For any undirected graph $\mathcal{G}(\mathcal{V},\mathcal{E})$, the cones $\mathbb{S}^n_{+}(\mathcal{E},?)$ and $\mathbb{S}_{+}^n(\mathcal{E},0)$ are dual to each other with respect to the trace inner product in the space of sparse matrices $\mathbb{S}^n(\mathcal{E},0)$~\cite{vandenberghe2014chordal}. In other words,
\begin{align*}
\mathbb{S}^n_{+}(\mathcal{E},?) &\equiv \{X\in\mathbb{S}^n(\mathcal{E},0):\; \langle X,Z\rangle\geq 0, \;\forall Z\in \mathbb{S}_{+}^n(\mathcal{E},0) \},
\\
\mathbb{S}^n_{+}(\mathcal{E},0) &\equiv \{Z\in\mathbb{S}^n(\mathcal{E},0):\; \langle Z,X\rangle\geq 0, \; \forall X\in \mathbb{S}_{+}^n(\mathcal{E},?) \}.
\end{align*}

If $\mathcal{G}$ is chordal, then $\mathbb{S}^n_{+}(\mathcal{E},?)$ and $\mathbb{S}_{+}^n(\mathcal{E},0)$ can be equivalently decomposed into a set of smaller but coupled convex cones according to the following theorems.
\begin{theorem} [{\!\cite[theorem 7]{grone1984positive}}]\label{T:ChordalCompletionTheorem}
     Let $\mathcal{G}(\mathcal{V},\mathcal{E})$ be a chordal graph and let $\{\mathcal{C}_1,\mathcal{C}_2, \ldots, \mathcal{C}_p\}$ be the set of its maximal cliques. Then, $X\in\mathbb{S}^n_+(\mathcal{E},?)$ if and only if
    $$ E_{\mathcal{C}_k} X E_{\mathcal{C}_k}^T \in \mathbb{S}^{\vert \mathcal{C}_k \vert}_+,
    \qquad k=1,\,\ldots,\,p.$$
\end{theorem}
\begin{theorem} [{\!\cite[theorem 2.3]{agler1988positive},~\cite[theorem 4]{griewank1984existence},~\cite[theorem 1]{kakimura2010direct}}]
\label{T:ChordalDecompositionTheorem}
     Let $\mathcal{G}(\mathcal{V},\mathcal{E})$ be a chordal graph and let $\{\mathcal{C}_1,\mathcal{C}_2, \ldots, \mathcal{C}_p\}$ be the set of its maximal cliques. Then, $Z\in\mathbb{S}^n_+(\mathcal{E},0)$ if and only if there exist matrices $Z_k \in \mathbb{S}^{\vert \mathcal{C}_k \vert}_+$ for $k=1,\,\ldots,\,p$ such that
    $$Z = \sum_{k=1}^{p} E_{\mathcal{C}_k}^T Z_k E_{\mathcal{C}_k}.$$
\end{theorem}
Note that these results can be proven individually, but can also be derived from each other using the duality of the cones $\mathbb{S}^n_{+}(\mathcal{E},?)$ and $\mathbb{S}_{+}^n(\mathcal{E},0)$~\cite{kim2011exploiting}.
In this paper, the terminology \emph{chordal} (or \emph{clique}) \emph{decomposition of a sparse matrix cone}  will refer to the application of Theorem~\ref{T:ChordalCompletionTheorem} or Theorem~\ref{T:ChordalDecompositionTheorem} to replace a large sparse PSD cone with a set of smaller but coupled PSD cones. Chordal decomposition of sparse matrix cones underpins much of the recent research on sparse SDPs~\cite{fukuda2001exploiting,kim2011exploiting,andersen2010implementation, sun2014decomposition,Madani2015ADMM,vandenberghe2014chordal}, most of which relies on the conversion framework for IPMs proposed in~\cite{fukuda2001exploiting, kim2011exploiting}.

To illustrate the concept, consider the chordal graph in Fig.~\ref{F:ChordalGraph}(b). {By} Theorem~\ref{T:ChordalCompletionTheorem},
$$
    \begin{bmatrix} X_{11} & X_{12} & 0 & X_{14} \\ X_{12} & X_{22} & X_{23} & X_{24}\\ 0 & X_{23} & X_{33} & X_{34} \\
    X_{14} & X_{24} & X_{34} & X_{44}\end{bmatrix} \in \mathbb{S}_+^n(\mathcal{E},?)
    \, \Leftrightarrow \,
    \begin{bmatrix}
    X_{11} & X_{12} & X_{14} \\ X_{12} & X_{22} & X_{24} \\
    X_{14} & X_{24} & X_{44}
    \end{bmatrix} \succeq 0,
    \quad
    \begin{bmatrix}
    X_{22} & X_{23} & X_{24} \\ X_{23} & X_{33} & X_{34} \\
    X_{24} & X_{34} & X_{44}
    \end{bmatrix} \succeq 0.
$$
Similarly, Theorem~\ref{T:ChordalDecompositionTheorem} guarantees that (after eliminating some of the variables)
$$
    \begin{bmatrix} Z_{11} & Z_{12} & 0 & Z_{14} \\ Z_{12} & Z_{22} & Z_{23} & Z_{24}\\ 0 & Z_{23} & Z_{33} & Z_{34} \\
    Z_{14} & Z_{24} & Z_{34} & Z_{44}\end{bmatrix} \in \mathbb{S}_+^n(\mathcal{E},0)
    \,\Leftrightarrow \,
    \begin{cases}
    \begin{bmatrix} Z_{11} & Z_{12} & Z_{14} \\ Z_{12} & {a_1} & a_3 \\
    Z_{14} & a_3 & {a_2}\end{bmatrix} \succeq 0,
    \quad
    \begin{bmatrix} {b_1} & Z_{23} & b_3 \\ Z_{23} & Z_{33} & Z_{34} \\
    b_3 & Z_{34} & {b_2}\end{bmatrix} \succeq 0,
    \\[2em]
    a_i + b_i = Z_{ii}, \, i \in {\{1,2\}},
    \\
    a_3 + b_3 = Z_{24}
    \end{cases}
$$
{for some constants $a_1$, $a_2$, $a_3$ and $b_1$, $b_2$, $b_3$.}
Note that the PSD contraints obtained after the chordal decomposition of $X$ (resp. $Z$) are coupled via the elements $X_{22}$, $X_{44},$ and $X_{24}=X_{42}$ (resp. $Z_{22}$, $Z_{44},$ and $Z_{24}=Z_{42}$).

\subsection{The Alternating Direction Method of Multipliers}

The computational ``engine'' employed in this work is the alternating direction method of multipliers (ADMM). ADMM is an operator-splitting method developed in the 1970s, and it is known to be equivalent to other operator-splitting methods such as Douglas-Rachford splitting and Spingarn's method of partial inverses; see~\cite{boyd2011distributed} for a review.
The ADMM algorithm solves the optimization problem
\begin{equation}
\label{e:ADMM_base_prob}
    \begin{aligned}
        \min_{x,y} \quad & f(x)+g(y) \\
        \text{subject to} \quad & Ax + By = c,
    \end{aligned}
\end{equation}
where $f$ and $g$ are convex functions, $x \in \mathbb{R}^{n_x}, y \in \mathbb{R}^{n_y}, A \in \mathbb{R}^{n_c\times n_x}, B \in \mathbb{R}^{n_c\times n_y}$ and $c \in \mathbb{R}^{n_c}$. Given a penalty parameter $\rho>0$ and a dual multiplier $z \in \mathbb{R}^{n_c}$, the ADMM algorithm finds a saddle point of the augmented Lagrangian
\begin{equation*}
 \mathcal{L}_{\rho}(x,y,z) := f(x) + g(y)
 + z^T \left( Ax + By - c  \right)
 + \frac{\rho}{2} \left\|Ax + By - c\right\|^2
\end{equation*}
by minimizing $\mathcal{L}$ with respect to the primal variables $x$ and $y$ separately, followed by a dual variable update:
\begin{subequations}\label{E:ADMM}
    \begin{align}
        x^{(n+1)} & = \text{arg} \min_{x} \; \mathcal{L}_{\rho}(x,y^{(n)},z^{(n)}),
        \label{E:ADMM_S1}\\
        y^{(n+1)} & = \text{arg} \min_{y} \; \mathcal{L}_{\rho}(x^{(n+1)},y,z^{(n)}),
        \label{E:ADMM_S2}\\
        z^{(n+1)} &= z^{(n)} + \rho\,( A x^{(n+1)} + B y^{(n+1)} - c). \label{E:ADMM_S3}\
    \end{align}
\end{subequations}
The superscript $(n)$ indicates that a variable is fixed to its value at the $n$-th iteration. Note that since $z$ is fixed in~\eqref{E:ADMM_S1} and~\eqref{E:ADMM_S2}, one may equivalently minimize the modified Lagrangian
\begin{equation*}
 \hat{\mathcal{L}}_{\rho}(x,y,z) := f(x) + g(y)
 + \frac{\rho}{2} \left\|Ax + By - c + \frac{1}{\rho} z\right\|^2.
\end{equation*}

Under very mild conditions, the ADMM converges to a solution of~\eqref{e:ADMM_base_prob} with a rate $\mathcal{O}(\frac{1}{n})$ ~\cite[Section 3.2]{boyd2011distributed}. ADMM is particularly suitable when~\eqref{E:ADMM_S1} and~\eqref{E:ADMM_S2} have closed-form expressions, or can be solved efficiently. Moreover, splitting the minimization over $x$ and $y$ often allows distributed and/or parallel implementations of steps \eqref{E:ADMM_S1}--\eqref{E:ADMM_S3}.

\section{Chordal decomposition of sparse SDPs} \label{sec:decompositionSDPs}

The sparsity pattern of the problem data for the primal-dual pair of standard-form SDPs~\eqref{E:PrimalSDP}-\eqref{E:DualSDP} can be described using the so-called \emph{aggregate sparsity pattern}. We say that the pair of SDPs~\eqref{E:PrimalSDP}-\eqref{E:DualSDP} has an aggregate sparsity pattern $\mathcal{G}(\mathcal{V},\mathcal{E})$ if
\begin{equation} \label{E:AggregateSparsity}
    C \in \mathbb{S}^n(\mathcal{E},0) \quad \text{and} \quad
    A_i \in \mathbb{S}^n(\mathcal{E},0), \; i = 1, \ldots, m.
\end{equation}
In other words, the aggregate sparsity pattern $\mathcal{G}$ is the union of the individual sparsity patterns of the data matrices $C$, $A_1,\,\ldots,\,A_m$. Throughout the rest of this paper, we assume that the aggregate sparsity pattern $\mathcal{G}$ is chordal (or that a suitable chordal extension has been found), and that it has $p$ maximal cliques $\mathcal{C}_1,\,\ldots,\,\mathcal{C}_p$. In addition, we assume that the matrices $A_1$, $\ldots$, $A_m$ are linearly independent.

It is not difficult to see that the aggregate sparsity pattern defines the sparsity pattern of any feasible dual variable $Z$ in~\eqref{E:DualSDP}, {\it i.e.}, any dual feasible $Z$ must have sparsity pattern $\mathcal{G}$. Similarly, while the primal variable $X$ in~\eqref{E:PrimalSDP} is usually dense, the value of the cost function and the equality constraints depend only on the entries $X_{ij}$ with $(i,j)\in\mathcal{E}$, and the remaining entries simply guarantee that {$X$ is PSD}. Recalling the definition of the sparse matrix cones $\mathbb{S}^n_{+}(\mathcal{E},?)$ and $\mathbb{S}^n_{+}(\mathcal{E},0)$, we can therefore recast the primal-form SDP~\eqref{E:PrimalSDP} as
\begin{equation}\label{E:NonSymmetricPrimal}
    \begin{aligned}
            \min_{X} \quad & \langle C, X \rangle \\
            \text{subject to} \quad & \langle A_i, X\rangle = b_i,
            \quad i = 1,\,\ldots,\,m, \\
                & X \in \mathbb{S}^n_{+}(\mathcal{E},?),
            \end{aligned}
\end{equation}
and the dual-form SDP~\eqref{E:DualSDP} as
\begin{equation}\label{E:NonSymmetricDual}
  \begin{aligned}
            \max_{y, Z} \quad & \langle b,y\rangle \\
            \text{subject to} \quad  & Z + \sum_{i=1}^m A_i\,y_i = C,\\
            & Z \in \mathbb{S}^n_{+}(\mathcal{E},0).
        \end{aligned}
\end{equation}
{This formulation} was first proposed by Fukuda \emph{et al.}~\cite{fukuda2001exploiting}, and was later discussed in~\cite{andersen2010implementation,sun2014decomposition, kim2011exploiting}. Note that~\eqref{E:NonSymmetricPrimal} and~\eqref{E:NonSymmetricDual} are a primal-dual pair of linear conic problems because the cones $\mathbb{S}^n_+(\mathcal{E}, ?)$ and $\mathbb{S}^n_+(\mathcal{E}, 0)$ are dual to each other.

\subsection{Domain-space decomposition}
\label{SS:DomainSpaceDecomp}

As we have seen in Section~\ref{sec:preliminaries}, Theorem~\ref{T:ChordalCompletionTheorem} allows us to decompose the sparse matrix cone constraint $X \in \mathbb{S}^n_{+}(\mathcal{E},?)$ into $p$ standard PSD constraints on the submatrices of $X$ defined by the cliques $\mathcal{C}_1,\,\ldots,\,\mathcal{C}_p$. In other words,
\begin{equation*}
    X \in \mathbb{S}^n_{+}(\mathcal{E},?)
    \, \Leftrightarrow \,
    E_{\mathcal{C}_k}XE_{\mathcal{C}_k}^T
    \in \mathbb{S}^{|\mathcal{C}_k|}_+,
    \quad k = 1,\, \ldots,\, p.
\end{equation*}
These $p$ constraints are implicitly coupled since $E_{\mathcal{C}_l}XE_{\mathcal{C}_l}^T$ and $E_{\mathcal{C}_q}XE_{\mathcal{C}_q}^T$ have overlapping elements if $\mathcal{C}_l \cap \mathcal{C}_q \neq \emptyset$. Upon introducing slack variables $X_k$, $k=1,\,\ldots,\,p$, we can rewrite this as
\begin{equation} \label{E:PrimalSlackVariables}
    X \in \mathbb{S}^n_{+}(\mathcal{E},?)
    \, \Leftrightarrow \,
    \begin{cases}
    		X_k = E_{\mathcal{C}_k}XE_{\mathcal{C}_k}^T, &k = 1,\,\ldots,\,p,\\
    		X_k \in \mathbb{S}^{|\mathcal{C}_k|}_+, &k = 1,\,\ldots,\,p.
    \end{cases}
\end{equation}
The primal optimization problem~\eqref{E:NonSymmetricPrimal} is then equivalent to the SDP
\begin{equation} \label{E:DecomposedPrimalSDP}
    \begin{aligned}
        \min_{X,X_1,\ldots,X_p} \quad & \langle C, X \rangle \\
        \text{subject to} \quad & \langle A_i, X\rangle = b_i, &i = 1, \ldots, m, \\
            & X_k = E_{\mathcal{C}_k}XE_{\mathcal{C}_k}^T,  &k = 1, \ldots, p, \\
            & X_k \in \mathbb{S}^{|\mathcal{C}_k|}_+,  &k = 1, \ldots, p.
     \end{aligned}
\end{equation}
Adopting the same terminology used in~\cite{fukuda2001exploiting}, we refer to~\eqref{E:DecomposedPrimalSDP} as the \emph{domain-space} decomposition of the primal-standard-form SDP~\eqref{E:PrimalSDP}.

\begin{remark}
\label{R:IPMvsFOM}
The main difference between the conversion method proposed in this section and that in~\cite{fukuda2001exploiting,kim2011exploiting} is that the large matrix $X$ is not eliminated. Instead, in the domain-space decomposition of~\cite{fukuda2001exploiting,kim2011exploiting}, $X$ is eliminated by replacing the constraints
\begin{equation*}
  X_k = E_{\mathcal{C}_k} X E_{\mathcal{C}_k}^T , \quad k=1,\,\ldots,\,p,
\end{equation*}
with the requirement that the entries of any two different sub-matrices $X_j,\,X_k$ must match if they map to the same entry in $X$. Mathematically, this condition can be written as
\begin{equation}\label{E:OverlappingEquality}
            E_{\mathcal{C}_j \cap \mathcal{C}_k}
            \left(
            		E_{\mathcal{C}_k}^T X_k E_{\mathcal{C}_k}
            	  - E_{\mathcal{C}_j}^T X_j E_{\mathcal{C}_j}
            	  \right)
            	 E_{\mathcal{C}_j \cap \mathcal{C}_k}^T =0,
            \quad \forall j,k \;\text{ such that }\; \mathcal{C}_j \cap \mathcal{C}_k \neq \emptyset.
\end{equation}
Redundant constraints in~\eqref{E:OverlappingEquality} can be eliminated using the \emph{running intersection property} of the cliques~\cite{blair1993introduction,fukuda2001exploiting}, and the decomposed SDP can be solved efficiently by IPMs in certain cases~\cite{fukuda2001exploiting,kim2011exploiting}. However, {applying FOMs to~\eqref{E:DecomposedPrimalSDP} effectively} after {the elimination of} $X$ is not straightforward {because the PSD matrix variables $X_1,\,\ldots,\,X_p$ are coupled via~\eqref{E:OverlappingEquality}}. In~\cite{sun2014decomposition}, {for example}, an SDP with a quadratic objective had to be solved at each iteration to impose the PSD constraints, requiring an additional iterative solver. Even when this problem is resolved, \emph{e.g.}, by using the algorithm of~\cite{ODonoghue2016}, the size of the KKT system enforcing the affine constraints is increased dramatically by the consensus conditions~\eqref{E:OverlappingEquality}, sometimes so much that memory requirements are prohibitive on desktop computing platforms~\cite{fukuda2001exploiting}.
In contrast, we show in Section \ref{sec:primal&dual} that if
a set of slack  variables $X_k$ are introduced in~\eqref{E:PrimalSlackVariables} and
{$X$ is not eliminated from~\eqref{E:DecomposedPrimalSDP}}, then the PSD constraint can be imposed via projections onto small PSD cones. At the same time, the affine constraints require the solution of an $m\times m$ linear system of equations, as if no consensus constraints were introduced. This makes our conversion framework more suitable for FOMs than that of~\cite{fukuda2001exploiting,kim2011exploiting},
{ as all steps in many common operator-splitting algorithms have an efficiently computable explicit solution. Of course, the equalities $X_k = E_{\mathcal{C}_k}XE_{\mathcal{C}_k}^T$, $k = 1, \ldots, p$ are satisfied only within moderate tolerances when FOMs are utilized, and the accumulation of small errors might make it more difficult to solve the original SDP to a given degree of accuracy compared to the methods in~\cite{fukuda2001exploiting, kim2011exploiting, sun2014decomposition, wen2010alternating}. Therefore, the trade-off between the gains in computational complexity and the reduction in accuracy should be carefully considered when choosing the most suitable approach to solve a given large-scale SDP. Nonetheless, our numerical experiments of Section~\ref{sec:simulation} demonstrate that working with~\eqref{E:DecomposedPrimalSDP} is often a 
 competitive strategy.}
\end{remark}

\subsection{Range-space decomposition}
\label{SS:RangeSpaceDecomp}

A \emph{range-space} decomposition of the dual-standard-form SDP~\eqref{E:DualSDP} can be formulated by applying Theorem~\ref{T:ChordalDecompositionTheorem} to the sparse matrix cone constraint  $Z\in \mathbb{S}^n_{+}(\mathcal{E},0)$ in~\eqref{E:NonSymmetricDual}:
\begin{equation*}
    Z \in \mathbb{S}^n_{+}(\mathcal{E},0)
    \, \Leftrightarrow \,
    Z = \sum_{k=1}^{p} E_{\mathcal{C}_k}^T Z_k E_{\mathcal{C}_k}
    {
    \text{ for some }
    Z_k \in \mathbb{S}^{|\mathcal{C}_k|}_+,\;
    k = 1,\,\ldots,\,p.
	}
\end{equation*}
We then introduce slack variables $V_k$, $k=1,\,\ldots,\,p$ and conclude that $Z\in \mathbb{S}^n_{+}(\mathcal{E},0)$ if and only if there exists matrices $Z_k,V_k \in \mathbb{S}^{|\mathcal{C}_k|}$, $k=1,\,\ldots,\,p$, such that
$$
Z = \sum_{k=1}^{p} E_{\mathcal{C}_k}^T V_k E_{\mathcal{C}_k}, \quad
Z_k = V_k, \;  k = 1,\,\ldots, \,p, \quad
Z_k \in \mathbb{S}^{|\mathcal{C}_k|}_+,  k = 1,\,\ldots, \,p.
$$
The range-space decomposition of~\eqref{E:DualSDP} is then given by
\begin{equation} \label{E:DecomposedDualSDP}
    \begin{aligned}
            \max_{y,Z_1,\ldots, Z_p, V_1,\ldots,V_p} \quad & \langle b,y\rangle \\
            \text{subject to} \quad  & \sum_{i=1}^m A_i\,y_i + \sum_{k=1}^{p} E_{\mathcal{C}_k}^T V_k E_{\mathcal{C}_k}= C,\\
            & Z_k - V_k = 0, \; k = 1, \ldots, p, \\
            & Z_k \in  \mathbb{S}^{|\mathcal{C}_k|}_+, \quad k = 1, \ldots, p.
        \end{aligned}
\end{equation}
{
Similar comments as in Remark~\ref{R:IPMvsFOM} hold:  the slack variables $V_1,\,\ldots,\,V_p$ are essential to formulate a decomposition framework suitable for the application of FOMs, although their introduction might complicate solving~\eqref{E:DualSDP} to a desired accuracy.}

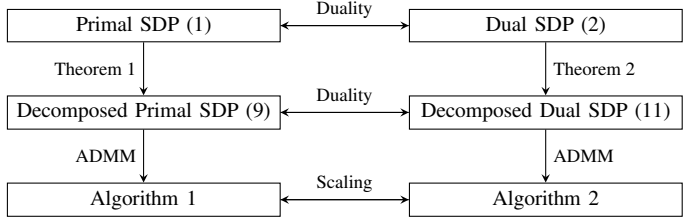
\begin{figure}
    \centering
    \setlength{\abovecaptionskip}{0pt}
    \setlength{\belowcaptionskip}{1em}
	\footnotesize
	\begin{tikzpicture}
	  \matrix (m) [matrix of nodes,
	  		       row sep = 2.5em,	
	  		       column sep = 6em,	
  			       nodes={rectangle, draw=black, align=center, text width=3.4cm}]
  	{
   	Primal SDP~\eqref{E:PrimalSDP} & Dual SDP~\eqref{E:DualSDP} \\
   	Decomposed Primal SDP~\eqref{E:DecomposedPrimalSDP} &
   	Decomposed Dual SDP~\eqref{E:DecomposedDualSDP} \\
   	Algorithm~\ref{A:ADMMPrimal} & Algorithm~\ref{A:ADMMDual}\\};
	\path[-stealth]
		(m-1-1) edge node [left, align=center]
            {\scriptsize Theorem~\ref{T:ChordalCompletionTheorem}}  (m-2-1)
		(m-1-2) edge node [right, align=center]
			{\scriptsize Theorem~\ref{T:ChordalDecompositionTheorem}}  (m-2-2)
		(m-1-1) edge node [above] {\scriptsize Duality} (m-1-2)
        (m-1-2) edge (m-1-1)
		(m-2-1) edge node [above] {\scriptsize Duality} (m-2-2)
        (m-2-2) edge (m-2-1)
		(m-2-1) edge node [left]  {\scriptsize ADMM}	   (m-3-1)
		(m-2-2) edge node [right] {\scriptsize ADMM}	   (m-3-2)
		(m-3-1) edge node [above] {\scriptsize Scaling} (m-3-2)
        (m-3-2) edge   (m-3-1);
	\end{tikzpicture}
    \caption{Duality between the original primal and dual SDPs, and the decomposed primal and dual SDPs.}
    \label{F:Duality}
\end{figure}

\begin{remark}
Although the domain- and range-space decompositions \eqref{E:DecomposedPrimalSDP} and \eqref{E:DecomposedDualSDP} have been derived individually, they are in fact a primal-dual pair of SDPs.  The duality between the original SDPs~\eqref{E:PrimalSDP} and~\eqref{E:DualSDP} is inherited by the decomposed SDPs~\eqref{E:DecomposedPrimalSDP} and~\eqref{E:DecomposedDualSDP} by virtue of the duality between Theorem~\ref{T:ChordalCompletionTheorem} and Theorem~\ref{T:ChordalDecompositionTheorem}. This elegant picture is illustrated in Fig.~\ref{F:Duality}.
\end{remark}

\section{ADMM for domain- and range-space decompositions of sparse SDPs}
\label{sec:primal&dual}

In this section, we demonstrate how ADMM can be applied to solve the domain-space decomposition~\eqref{E:DecomposedPrimalSDP} and the range-space decomposition~\eqref{E:DecomposedDualSDP} efficiently. Furthermore, we show that the resulting domain- and range-space algorithms are equivalent, in the sense that one is just a scaled version of the other {(cf. Fig.~\ref{F:Duality})}.
Throughout this section, $\delta_{\mathcal{K}}(x)$ will denote the indicator function of a set $\mathcal{K}$, {\it i.e.},
$$
    \delta_{\mathcal{K}}(x) =
    \begin{cases}
    0, &\text{if } x \in \mathcal{K}, \\
    + \infty, &\text{otherwise}.
    \end{cases}
$$
For notational neatness, however, we write $\delta_{0}$ when $\mathcal{K}\equiv \{0\}$.

To ease the exposition further, we consider the usual vectorized forms of~\eqref{E:DecomposedPrimalSDP} and~\eqref{E:DecomposedDualSDP}. Specifically, we let $\vect:\mathbb{S}^n \to \mathbb{R}^{n^2}$ be the usual operator mapping a matrix to the stack of its columns and define the vectorized data
\begin{align*}
c := \vect(C),   \quad
A := \begin{bmatrix} \vect(A_0) & \hdots & \vect(A_m) \end{bmatrix}^T.
\end{align*}
Note that the assumption that $A_1$, $\ldots$, $A_m$ are linearly independent matrices means that $A$ has  full row rank. For all $k = 1,\,\ldots,\,p$, we also introduce the vectorized variables
\begin{align*}
        x := \vect(X), \quad
        x_k := \vect(X_k), \quad
        z_k := \vect(Z_k),  \quad
        v_k := \vect(V_k),
\end{align*}
and define ``entry-selector'' matrices $H_k := E_{\mathcal{C}_k} \otimes E_{\mathcal{C}_k}$ for $k = 1, \ldots, p$
that project $x$ onto the subvectors $x_1,\,\ldots,\,x_p$, {\it i.e.}, such that
\begin{equation*}
    x_k = \vect(X_k) = \vect(E_{\mathcal{C}_k} X E_{\mathcal{C}_k}^T) = H_k x.
\end{equation*}
Note that for each $k=1,\,\ldots,\,p$, the rows of $H_k$ are orthonormal, and that the matrix $H_k^TH_k$ is diagonal. Upon defining
\begin{equation*}
\mathcal{S}_k := \left\{
x \in \mathbb{R}^{\vert \mathcal{C}_k \vert^2}:\;
\vect^{-1}(x) \in \mathbb{S}_+^{\vert \mathcal{C}_k \vert}
\right\},
\end{equation*}
such that $x_k\in\mathcal{S}_k$ if and only if $X_k\in\mathbb{S}^{|\mathcal{C}_k|}_+$, we can rewrite~\eqref{E:DecomposedPrimalSDP} as
\begin{equation}\label{E:DecomposedPrimalVector}
        \begin{aligned}
        \min_{x,x_1,\ldots,x_p} \quad & \langle c,x \rangle\\
        \text{subject to} \quad & Ax=b,\\
        					   & x_k = H_k x, &k=1,\,\ldots,\,p,\\
        					   & x_k \in \mathcal{S}_k, &k=1,\,\ldots,\,p,
        \end{aligned}
\end{equation}
while~\eqref{E:DecomposedDualSDP} becomes
\begin{equation}\label{E:DecomposedDualVector}
  \begin{aligned}
            \max_{y,z_1,\ldots,z_p,v_1,\ldots,v_p} \quad & \langle b,y \rangle\\
        \text{subject to} \quad
        & A^Ty + \sum_{k=1}^p H_k^Tv_k =  c ,\\
        & z_k - v_k = 0, & k=1,\,\ldots,\,p,\\
        & z_k \in \mathcal{S}_k, & k=1,\,\ldots,\,p.
        \end{aligned}
\end{equation}

\subsection{ADMM for the domain-space decomposition} \label{subse:primalADMM}

We start by moving the constraints $Ax = b$ and $x_k \in \mathcal{S}_k$ in \eqref{E:DecomposedPrimalVector} to the objective using the indicator functions $\delta_0(\cdot)$ and $\delta_{\mathcal{S}_k}(\cdot)$, respectively, \textit{i.e.}, we write
\begin{equation} \label{E:ADMMPrimal}
    \begin{aligned}
    \min_{x,x_1,\ldots,x_p} \quad &\langle c,x \rangle + \delta_0\left( Ax - b \right) + \sum_{k=1}^{p} \delta_{\mathcal{S}_k}(x_k)
    \\
    \text{subject to} \quad &x_k = H_k x, \quad k=1,\,\ldots,\,p.
    \end{aligned}
\end{equation}

This problem is in the standard form for the application of ADMM. Given a penalty parameter $\rho>0$ and a Lagrange multiplier $\lambda_k$ for each constraint $x_k = H_k x$, $k=1,\,\ldots,\,p$, we consider the (modified) augmented Lagrangian
\begin{multline}
\label{E:AugLagrPrimal}
\mathcal{L}(x,x_1,\ldots,x_k,\lambda_1,\ldots,\lambda_k):= \langle c,x \rangle + \delta_0\left( Ax - b \right)
\\
+ \sum_{k=1}^{p} \left[ \delta_{\mathcal{S}_k}(x_k) + \frac{\rho}{2}\left\| x_k - H_k x + \frac{1}{\rho}\lambda_k \right\|^2 \right],
\end{multline}
and group the variables as $\mathcal{X}:= \{x\}$,
$\mathcal{Y}:= \{x_1,\,\ldots,\,x_p\}$, and
$\mathcal{Z}:= \{\lambda_1,\,\ldots,\,\lambda_p\}$.
According to \eqref{E:ADMM}, each iteration of the ADMM requires the minimization of the Lagrangian in~\eqref{E:AugLagrPrimal} with respect to the $\mathcal{X}$- and $\mathcal{Y}$-blocks separately, {followed} by an update of the multipliers $\mathcal{Z}$. At each step, the variables not being optimized over are fixed to their most current value. Note that splitting the primal variables $x,\,x_1, \ldots ,x_p$ in the two blocks $\mathcal{X}$ and $\mathcal{Y}$ defined above is essential to solving the $\mathcal{X}$ and $\mathcal{Y}$ minimization sub-problems~\eqref{E:ADMM_S1} and~\eqref{E:ADMM_S2}; more details will be given in Remark~\ref{Remark:PrimalSplitting} after describing the $\mathcal{Y}$-minimization step in Section~\ref{se:MinXblkprimal}.

\subsubsection{Minimization over $\mathcal{X}$} \label{se:MinYblkprimal}

Minimizing the augmented Lagrangian~\eqref{E:AugLagrPrimal} over $\mathcal{X}$ is equivalent to the equality-constrained quadratic program
\begin{equation} \label{E:MinXblockPrimal}
    \begin{aligned}
        \min_{x} \quad &\langle c,x \rangle + \frac{\rho}{2}\sum_{k=1}^{p} \left\| x_k^{(n)} - H_k x + \frac{1}{\rho}\lambda_k^{(n)} \right\|^2
        \\
        \text{subject to} \quad & Ax=b.
    \end{aligned}
\end{equation}
Letting $\rho y$ be the multiplier for the equality constraint (we scale the multiplier by $\rho$ for convenience), and defining
\begin{equation}
D := \sum_{k=1}^{p} H_k^T H_k,
\end{equation}
the optimality conditions for~\eqref{E:MinXblockPrimal} can be written as the KKT system
\begin{equation} \label{E:OptCondMinYPrimal}
    \begin{bmatrix}D & A^T \\ A & 0\end{bmatrix}
    \begin{bmatrix}x \\ y\end{bmatrix} =
    \begin{bmatrix}\sum_{k=1}^{p} H_k^T\left( x_k^{(n)}+\rho^{-1}\lambda_k^{(n)}\right) - \rho^{-1}c \\ b \end{bmatrix}.
\end{equation}
Recalling that the product $H_k^T H_k$ is a diagonal matrix for all $k=1,\,\ldots,\,p$  we conclude that so is $D$, and since $A$ has full row rank by assumption~\eqref{E:OptCondMinYPrimal} can be solved efficiently, for instance by block elimination. In particular, eliminating $x$ shows that the only matrix to be inverted/factorized is
\begin{equation} \label{E:ADMMprimalFatorization}
    AD^{-1}A^T \in \mathbb{S}^{m}.
\end{equation}
Incidentally, we note that the first-order algorithms of~\cite{wen2010alternating,ODonoghue2016} require the factorization of a similar matrix with the same dimension.
Since this matrix is the same at every iteration, its Cholesky factorization (or any other factorization of choice) can be computed and cached before starting the ADMM iterations. For some families of SDPs, such as the SDP relaxation of MaxCut problems and sum-of-squares (SOS) feasibility problems~\cite{zheng2017Exploiting}, the matrix $AD^{-1}A^T$ is diagonal, so solving~\eqref{E:OptCondMinYPrimal} is inexpensive even when the SDPs are very large. If factorizing $AD^{-1}A^T$ is too expensive, the linear system~\eqref{E:OptCondMinYPrimal} can alternatively be solved by an iterative method, such as the conjugate gradient method~\cite{saad2003iterative}.

\subsubsection{Minimization over $\mathcal{Y}$} \label{se:MinXblkprimal}
Minimizing the augmented Lagrangian~\eqref{E:AugLagrPrimal} over $\mathcal{Y}$ is equivalent to solving $p$ independent conic problems of the form
\begin{equation}
    \begin{aligned}
    \min_{x_k} \quad &\left\| x_k - H_k x^{(n+1)}  + {\rho}^{-1}\lambda_k^{(n)} \right\|^2
    \\
    \text{subject to} \quad &x_k \in \mathcal{S}_k.
    \end{aligned}
\end{equation}
In terms of the original matrix variables $X_1,\,\ldots,\,X_p$, each of these $p$ sub-problems amounts to a projection on a PSD cone. More precisely, if $\mathbb{P}_{\mathbb{S}^{\vert \mathcal{C}_k\vert}_+}$ denotes the projection onto the PSD cone $\mathbb{S}^{\vert \mathcal{C}_k\vert}_{+}$ and $\mat(\cdot) =\vect^{-1}(\cdot)$, we have
\begin{equation} \label{E:XkUpdate}
    x_k^{(n+1)} = \vect\left\{
    \mathbb{P}_{\mathbb{S}^{\vert \mathcal{C}_k\vert}_{+}} \left[
    \mat\left( H_k x^{(n+1)} -  {\rho}^{-1}\lambda_k^{(n)}\right)
    \right] \right\}
    .
\end{equation}

Since the size of each cone $\mathbb{S}^{\vert \mathcal{C}_k\vert}_{+}$ is small for typical sparse SDPs
and the projection onto it can be computed with an eigenvalue decomposition, the variables $x_1,\,\ldots,\,x_p$ can be updated efficiently. Moreover, the computation can be carried out in parallel. In contrast, the algorithms for generic SDPs developed in~\cite{wen2010alternating,ODonoghue2016,malick2009regularization} require projections onto the {(much larger)} original PSD cone $\mathbb{S}^n_+$.

\begin{remark} \label{Remark:PrimalSplitting}
    As anticipated in Remark~\ref{R:IPMvsFOM}, retaining the global variable $x$ in the domain-space decomposed SDP to enforce the consensus constraints between the entries of the subvectors $x_1,\ldots,x_p$ (\emph{i.e.}, $x_k = H_kx$) is fundamental. In fact, it allowed us to separate the conic constraints from the affine constraints in~\eqref{E:DecomposedPrimalVector} when applying the splitting strategy of ADMM, making the minimization over $\mathcal{Y}$ easy to compute and parallelizable. In contrast, when $x$ is eliminated as in the conversion method of~\cite{fukuda2001exploiting,kim2011exploiting}, the conic constraints and the affine constraints cannot be easily decoupled when applying the first-order splitting method: in~\cite{sun2014decomposition} a quadratic SDP had to be solved at each iteration, which limits its scalability.
\end{remark}

\subsubsection{Updating the multipliers $\mathcal{Z}$}
The final step in the $n$-th ADMM iteration is to update the multipliers $\lambda_1,\,\ldots,\,\lambda_p$ with the usual gradient ascent rule: for each $k=1,\,\ldots,\,p$,
\begin{equation} \label{E:LambdakUpdate}
    \lambda_k^{(n+1)} =
    \lambda_k^{(n)} + \rho \left(x_k^{(n+1)} -  H_k x^{(n+1)} \right).
\end{equation}
This computation is {inexpensive} and easily parallelized.

\subsubsection{Stopping conditions} \label{se:PrimalStoppingCond}
The ADMM algorithm is stopped after the $n$-th iteration if the relative primal/dual error measures.
\begin{subequations}
\begin{align}
    {\epsilon_\mathrm{c}} &= \frac{\left( \displaystyle\sum_{k=1}^p \left\| x_k^{(n)} -  H_k x^{(n)} \right\|^2 \right)^{1/2}}
    {
    \max\left\{
    \left(\displaystyle\sum_{k=1}^p \left\| x_k^{(n)}\right\|^2\right)^{1/2},
     	\left(\displaystyle\sum_{k=1}^p \left\| H_k x^{(n)}\right\|^2\right)^{1/2}
    \right\}
    }
, \label{E:PrimalConsensusResidual}\\[1em]
 {\epsilon_\lambda} &= \rho\, \left( \displaystyle\sum_{k=1}^p \left\| x_k^{(n)} -  x_k^{(n-1)} \right\|^2 \right)^{1/2}\left(\displaystyle\sum_{k=1}^p \left\| \lambda_k^{(n)}\right\|^2\right)^{-1/2},
    \end{align}
\end{subequations}
are smaller than a specified tolerance, $\epsilon_\mathrm{tol}$. The reader is referred to~\cite{boyd2011distributed} for a detailed discussion of stopping conditions for ADMM algorithms. In conclusion, a primal-form SDP with domain-space decomposition \eqref{E:DecomposedPrimalVector} can be solved using the steps summarized in Algorithm~\ref{A:ADMMPrimal}.

\begin{algorithm}[t]
\caption{ADMM for the domain-space decomposition of sparse {primal-form} SDPs}
\label{A:ADMMPrimal}
\begin{algorithmic}[1]
\algsize
\State Set $\rho>0$, $\epsilon_\mathrm{tol} >0$, a maximum number of iterations $n_{\max}$, and initial guesses $x^{(0)}$, $x_1^{(0)},\,\ldots,\,x_p^{(0)}$, $\lambda_1^{(0)},\,\ldots,\,\lambda_p^{(0)}$.
\State Data preprocessing: chordal extension, chordal decomposition, and factorization of the KKT system \eqref{E:OptCondMinYPrimal}.
\For{$n = 1,2, \ldots, n_{\max}$}
    \State Compute $x^{(n)}$ using~\eqref{E:OptCondMinYPrimal}.
    {
	\For{$k=1,\,\ldots,\,p$}
		\State Compute $x_k^{(n)}$ using~\eqref{E:XkUpdate}.
		\State Compute $\lambda_k^{(n)}$ using~\eqref{E:LambdakUpdate}.
	\EndFor
	}
	\State Update the residuals $\epsilon_\mathrm{c}, \epsilon_\lambda$.
    \If{$\max(\epsilon_\mathrm{c}, \epsilon_\lambda) \leq \epsilon_\mathrm{tol}$}
            \State \textbf{break}
        \EndIf
\EndFor
\end{algorithmic}
\end{algorithm}

\subsection{ADMM for the range-space decomposition}

An ADMM algorithm similar to Algorithm~\ref{A:ADMMPrimal} can be developed for the range-space decomposition~\eqref{E:DecomposedDualVector} of a dual-standard-form sparse SDP. As in Section~\ref{subse:primalADMM}, we start by moving all but the consensus equality constraints $z_k=v_k$, $k=1,\,\ldots,\,p$, to the objective using indicator functions. This leads to
\begin{align}
\label{E:DualDecomposedSDPvec}
\min \quad
&-\langle b,y \rangle + \delta_0\left( c-A^T y-\sum_{k=1}^pH_k^T v_k\right)
+ \sum_{k=1}^{p} \delta_{\mathcal{S}_k}(z_k)
\notag \\
\text{subject to} \quad &
z_k = v_k, \quad k=1,\,\ldots,\,p.
\end{align}

Given a penalty parameter $\rho>0$ and a Lagrange multiplier $\lambda_k$ for each of the constraints $z_k = v_k$, $k=1,\,\ldots,\,p$, we consider the (modified) augmented Lagrangian
\begin{multline}
\mathcal{L}(y,v_1,\ldots,v_p,z_1,\ldots,z_p,\lambda_1,\ldots,\lambda_p):=
-\langle b, y \rangle
\\
+ \delta_0\left( c - A^T y - \sum_{k=1}^p H_k^T v_k \right)
+ \sum_{k=1}^{p} \left[ \delta_{\mathcal{S}_k}(z_k) + \frac{\rho}{2}\left\| z_k - v_k + \frac{1}{\rho}\lambda_k \right\|^2 \right],
\label{E:AugLagrDual}
\end{multline}
and consider three groups of variables,
$\mathcal{X}:= \{y,v_1,\,\ldots,\,v_p\}$,
$\mathcal{Y}:= \{z_1,\,\ldots,\,z_p\}$, and
$\mathcal{Z}:= \{\lambda_1,\,\ldots,\,\lambda_p\}$.
Similar to Section~\ref{subse:primalADMM}, each iteration of the ADMM algorithm for~\eqref{E:DecomposedDualVector} consists of minimizations over $\mathcal{X}$ and $\mathcal{Y}$, and an update of the  multipliers $\mathcal{Z}$. Each of these steps admits an inexpensive closed-form solution, as we demonstrate next.

\subsubsection{Minimization over $\mathcal{X}$}
Minimizing~\eqref{E:AugLagrDual} over block $\mathcal{X}$ is equivalent to solving the equality-constrained quadratic program
\begin{align}
\label{E:MinYblockDual}
\min_{y,v_1,\ldots,v_p} \quad &-\langle b,y\rangle + \frac{\rho}{2}\sum_{k=0}^{p} \left\| z_k^{(n)} - v_k + \frac{1}{\rho}\lambda_k^{(n)} \right\|^2
\notag\\
\text{subject to} \quad &c - A^T y - \sum_{k=1}^p H_k^Tv_k = 0.
\end{align}
Let $\rho x$ be the multiplier for the equality constraint. After some algebra, the optimality conditions for~\eqref{E:MinYblockDual} can be written as the KKT system
    \begin{equation}
    \label{E:OptCondMinYDual}
    \begin{bmatrix}D & A^T \\ A & 0\end{bmatrix}
    \begin{bmatrix}x \\ y\end{bmatrix} =
    \begin{bmatrix}c - \sum_{k=1}^{p} H_k^T\left( z_k^{(n)}+\rho^{-1}\lambda_k^{(n)}\right) \\ -\rho^{-1}b \end{bmatrix},
    \end{equation}
    plus a set of $p$ uncoupled equations for the variables $v_k$,
    \begin{equation}
    \label{E:vkEqn}
    v_k = z_k^{(n)} +\frac{1}{\rho} \lambda_k^{(n)} + H_k x,  \quad k=1,\,\ldots,\,p.
    \end{equation}

The KKT system~\eqref{E:OptCondMinYDual} is the same as~\eqref{E:OptCondMinYPrimal} after rescaling $x\mapsto-x$, $y\mapsto -y$, $c\mapsto \rho^{-1}c$ and $b\mapsto \rho b$. Consequently, the numerical cost {of~\eqref{E:MinYblockDual}} is the same as in Section~\ref{se:MinYblkprimal} plus the cost of~\eqref{E:vkEqn}, which is {inexpensive} and can be parallelized. Moreover, as in Section~\ref{se:MinYblkprimal}, the factors of the coefficient matrix required to solve the KKT system~\eqref{E:OptCondMinYDual} can be pre-computed and cached before iterating the ADMM algorithm.

\subsubsection{Minimization over $\mathcal{Y}$}
\label{se:MinXblk_dual}

As in Section~\ref{se:MinXblkprimal}, the variables $z_1,\,\ldots,\,z_p$ are updated with $p$ independent projections,
\begin{equation}
\label{E:zkUpdate}
z_k^{(n+1)} = \vect\left\{ \mathbb{P}_{\mathbb{S}^{\vert \mathcal{C}_k\vert}_+} \left[ \mat\left( v_k^{(n+1)} -  {\rho}^{-1}\lambda_k^{(n)}\right) \right] \right\},
\end{equation}
where $\mathbb{P}_{\mathbb{S}^{\vert \mathcal{C}_k\vert}_+} $ denotes projection on the PSD cone $\mathbb{S}^{\vert \mathcal{C}_k\vert}_{+}$. Again, these projections can be computed efficiently and in parallel.

\begin{remark}
As anticipated in Section~\ref{SS:RangeSpaceDecomp}, introducing the set of slack variables $v_k$ and the consensus constraints $z_k =v_k $, $k=1,\,\ldots,\,p$ is essential to obtain an efficient  algorithm for range-space decomposed SDPs. The reason is that the splitting strategy of the ADMM decouples the conic and affine constraints, and the conic variables can be updated using the simple conic projection~\eqref{E:zkUpdate}.
\end{remark}

\subsubsection{Updating the multipliers $\mathcal{Z}$}
The multipliers $\lambda_k$, $k=1,\,\ldots,\,p$, are updated (possibly in parallel) with the {inexpensive} gradient ascent rule
\begin{equation}
\label{E:MultUpdateDualSDP}
\lambda_k^{(n+1)} = \lambda_k^{(n)} + \rho\left( z_k^{(n+1)} - v_k^{(n+1)}\right).
\end{equation}

\subsubsection{Stopping conditions}
\label{se:DualStoppingCond}

Similar to Section~\ref{se:PrimalStoppingCond}, we stop our ADMM algorithm after the $n$-th iteration if the relative primal/dual error measures
\begin{subequations}
\begin{align}
{\epsilon_\mathrm{c}} &= \frac{\left( \displaystyle\sum_{k=1}^p \left\| z_k^{(n)} -  v_k^{(n)} \right\|^2 \right)^{1/2}}
{\max\left\{ \left(\displaystyle\sum_{k=1}^p \left\| z_k^{(n)}\right\|^2\right)^{1/2},
 		     \left(\displaystyle\sum_{k=1}^p \left\| v_k^{(n)}\right\|^2\right)^{1/2} \right\}},
\label{E:DualConsensusResidual}\\[1em]
{\epsilon_\lambda} &= \rho\,
\left(
\displaystyle\sum_{k=1}^p \left\| z_k^{(n)} -  z_k^{(n-1)} \right\|^2
\right)^{1/2}
\left(\displaystyle\sum_{k=1}^p \left\| \lambda_k^{(n)}\right\|^2\right)^{-1/2},
\end{align}
\end{subequations}
are smaller than a specified tolerance, $\epsilon_\mathrm{tol}$. The ADMM algorithm to solve the range-space decomposition~\eqref{E:DecomposedDualVector} of a dual-form sparse SDP is summarized in Algorithm~\ref{A:ADMMDual}.

\begin{algorithm}[t]
    \caption{{ADMM for the range-space decomposition of sparse dual-form SDPs}}
    \label{A:ADMMDual}
    \begin{algorithmic}[1]
    \algsize
    \State Set $\rho>0$, $\epsilon_\mathrm{tol} >0$, a maximum number of iterations $n_{\max}$ and initial guesses $y^{(0)}$, $z_1^{(0)},\,\ldots,\,z_p^{(0)}$, $\lambda_1^{(0)},\,\ldots,\,\lambda_p^{(0)}$.
    \State Data preprocessing: chordal extension, chordal decomposition, and factorization of the KKT system \eqref{E:OptCondMinYDual}.
    \For{$n = 1, 2, \ldots, n_{\max}$}
    	\For{$k=1,\,\ldots,\,p$}
    		\State Compute $z_k^{(n)}$ using~\eqref{E:zkUpdate}.
    	\EndFor
    	\State Compute $y^{(n)},x$ using~\eqref{E:MinYblockDual}.
    	\For{$k=1,\,\ldots,\,p$}
    		\State Compute $v_k^{(n)}$ using~\eqref{E:vkEqn}
            \vspace{0.3mm}
    		\State Compute $\lambda_k^{(n)}$ using~\eqref{E:DualMultUpdate} (no cost).
    	\EndFor
    	\State Update the residuals $\epsilon_\mathrm{c}$ and $\epsilon_\lambda$.
        \If{$\max(\epsilon_\mathrm{c}, \epsilon_\lambda) \leq \epsilon_\mathrm{tol}$}
            \State \textbf{break}
        \EndIf
    \EndFor
\end{algorithmic}
\end{algorithm}

\subsection{Equivalence between the primal and dual ADMM algorithms}

Since the computational cost of~\eqref{E:vkEqn} is the same as~\eqref{E:LambdakUpdate}, all ADMM iterations for the dual-form SDP with range-space decomposition~\eqref{E:DecomposedDualVector} have the same cost as those for the primal-form SDP with domain-space decomposition~\eqref{E:DecomposedPrimalVector}, plus the cost of~\eqref{E:MultUpdateDualSDP}.
However, if one minimizes the dual augmented Lagrangian~\eqref{E:AugLagrDual} over $z_1,\ldots,z_p$ \emph{before} minimizing it over $y,v_1,\ldots,v_p$, then~\eqref{E:vkEqn} can be used to simplify the multiplier update equations to
    \begin{equation} \label{E:DualMultUpdate}
        \lambda_k^{(n+1)} =  \rho H_k x^{(n+1)}, \quad k=1,\,\ldots,\,p.
    \end{equation}
Given that the products $H_1 x,\ldots,H_px$ have already been computed to update $v_1,\ldots,v_p$ in~\eqref{E:vkEqn}, updating the multipliers $\lambda_1,\ldots,\lambda_p$ requires only a scaling operation. {Then, after swapping the order of $\mathcal{X}$- and $\mathcal{Y}$-block minimization of~\eqref{E:AugLagrDual} and recalling that~\eqref{E:OptCondMinYPrimal} and~\eqref{E:OptCondMinYDual} are scaled versions of the same KKT system, the ADMM algorithms for the primal and dual standard form SDPs can be considered scaled versions of each other}; see Fig.~\ref{F:Duality} for an illustration. In fact, the equivalence between ADMM algorithms for the original ({\it i.e.}, before chordal decomposition) primal and dual SDPs was already noted in~\cite{yan2016self}.

\begin{remark}
    Although the iterates of Algorithm~\ref{A:ADMMPrimal} and Algorithm~\ref{A:ADMMDual} are the same up to scaling, the convergence performance of these two algorithms differ in practice because first-order methods are sensitive to the scaling of the problem data and of the iterates.
\end{remark}

\section{Homogeneous self-dual embedding of domain- and range-space decomposed SDPs}
\label{sec:hsde}

Algorithms~\ref{A:ADMMPrimal} and~\ref{A:ADMMDual}, as well as other first-order algorithms that exploit chordal sparsity~\cite{sun2014decomposition, Kalbat2015Fast, Madani2015ADMM}, can solve feasible problems, but cannot detect infeasibility in their current formulation. Although some recent ADMM methods resolve this issue~\cite{BGSBfeasibility2017,liu2017new}, an elegant way to deal with an infeasible primal-dual pair of  SDPs---which we pursue here---is to solve their homogeneous self-dual embedding (HSDE)~\cite{ye1994nl}.

The essence of the HSDE method is to search for a non-zero point in the intersection of a convex cone and a linear space; this is non-empty because it always contains the origin, meaning that the problem is always feasible. Given such a non-zero point, one can either recover optimal primal and dual solutions of the original pair of optimization problems, or construct a certificate of primal or dual infeasibility. HSDEs have been widely used to develop IPMs for SDPs~\cite{sturm1999using,ye2011interior}, and more recently O'Donoghue \emph{et al}. have proposed an operator-splitting method to solve the HSDE of general conic programs~\cite{ODonoghue2016}.

In this section, we formulate the HSDE of the domain- and range-space decomposed SDPs~\eqref{E:DecomposedPrimalVector} and~\eqref{E:DecomposedDualVector}, which is a primal-dual pair of SDPs.  We also apply ADMM to solve this HSDE; in particular, we extend the algorithm of~\cite{ODonoghue2016} to exploit chordal sparsity without increasing its computational cost (at least to leading order) compared to Algorithms~\ref{A:ADMMPrimal} and~\ref{A:ADMMDual}.

\subsection{Homogeneous self-dual embedding}
To simplify the formulation of the HSDE of the decomposed (vectorized) SDPs~\eqref{E:DecomposedPrimalVector} and~\eqref{E:DecomposedDualVector}, we let
$\mathcal{S} := \mathcal{S}_1 \times \cdots \times \mathcal{S}_p$
be the direct product of all semidefinite cones and define
\begin{align*}
  s := \begin{bmatrix}x_1 \\ \vdots\\ x_p \end{bmatrix}, \quad
  z := \begin{bmatrix}z_1 \\ \vdots\\ z_p \end{bmatrix}, \quad
  t := \begin{bmatrix}v_1 \\ \vdots\\ v_p \end{bmatrix}, \quad
  H := \begin{bmatrix}H_1 \\ \vdots\\ H_p \end{bmatrix}.
\end{align*}

When strong duality holds, the tuple $(x^*,s^*,y^*,t^*,z^*)$ is optimal if and only if all of the following conditions hold:
\begin{enumerate}
  \item  $(x^*,s^*)$ is primal feasible, {\it i.e.}, $Ax^*=b$, $s^*=Hx^*$, and $s^*\in\mathcal{S}$. For reasons that will become apparent below, we introduce slack variables $r^*=0$ and $w^*=0$ of appropriate dimensions and rewrite these conditions as
        \begin{equation}\label{E:PrimalFeasible}
            \begin{aligned}
                 Ax^* - r^* &= b, &
                 s^* + w^* &= Hx^*, &
                 s^* &\in \mathcal{S}, &
                 r^* &= 0,  &&&
                 w^* &= 0.
            \end{aligned}
        \end{equation}
  \item $(y^*,t^*,z^*)$ is dual feasible, {\it i.e.}, $A^Ty^* + H^Tt^* =  c$, $z^*=t^*$, and $z^*\in\mathcal{S}$. Again, it is convenient to introduce a slack variable $h^*=0$ of appropriate size and write
        \begin{equation}\label{E:DualFeasible}
            \begin{aligned}
                 A^Ty^* + H^Tt^* + h^* &=  c , &&&
                z^* - t^* &= 0, &&&
                z^* &\in \mathcal{S},&&&
                h^* &= 0.
            \end{aligned}
        \end{equation}
  \item The duality gap is zero,  \textit{i.e.}
        \begin{equation}\label{E:ZeroGap}
            c^Tx^* - b^Ty^* = 0.
        \end{equation}
\end{enumerate}

The idea behind the HSDE~\cite{ye1994nl} is to introduce two non-negative and complementary variables $\tau$ and $\kappa$ and embed the optimality conditions \eqref{E:PrimalFeasible}, \eqref{E:DualFeasible} and \eqref{E:ZeroGap} into the linear system $v = Q u$ with $u$, $v$ and $Q$ defined as
\begin{align} \label{E:HSDEnotations}
  u := \begin{bmatrix}
    x \\ s \\ y \\ t \\ \tau
  \end{bmatrix},
  \quad
  v := \begin{bmatrix}
    h \\ z \\ r \\ w \\ \kappa
  \end{bmatrix},
  \quad
  Q := \begin{bmatrix}
    0 & 0 & -A^T & -H^T & c \\
    0 & 0 & 0 & I & 0\\
    A & 0 & 0 & 0 & -b\\
    H & -I & 0 & 0 & 0 \\
    -c^T & 0  & b^T & 0 & 0\\
  \end{bmatrix}.
\end{align}
Any nonzero solution of this embedding can be used to recover an optimal solution for \eqref{E:DecomposedPrimalSDP} and \eqref{E:DecomposedDualSDP}, or provide a certificate for primal or dual infeasibility, depending on the values of $\tau$ and $\kappa$; details are omitted for brevity, and the interested reader is referred to~\cite{ODonoghue2016}.

The decomposed primal-dual pair of (vectorized) SDPs~\eqref{E:DecomposedPrimalVector}-\eqref{E:DecomposedDualVector} can therefore be recast as the self-dual conic feasibility problem
\begin{equation} \label{E:ADMMForm}
\begin{aligned}
\text{find} \quad &(u,v)\\
\text{subject to} \quad &v = Qu,\\
			      &(u,v) \in \mathcal{K} \times \mathcal{K}^*,
\end{aligned}
\end{equation}
where, writing $n_d =\sum_{k=1}^p |\mathcal{C}_k|^2$ for brevity,
$\mathcal{K} := \mathbb{R}^{n^2} \times \mathcal{S} \times \mathbb{R}^{m} \times \mathbb{R}^{n_d} \times \mathbb{R}_{+}$ is a cone and $\mathcal{K}^* := \{0\}^{n^2} \times \mathcal{S} \times \{0\}^{m} \times \{0\}^{n_d} \times \mathbb{R}_{+}$ is its dual.

\subsection{A simplified ADMM algorithm}

The feasibility problem~\eqref{E:ADMMForm} is in a form suitable for the application of ADMM, and moreover steps~\eqref{E:ADMM_S1}-\eqref{E:ADMM_S3} can be greatly simplified by virtue of its self-dual character~\cite{ODonoghue2016}. Specifically, the $n$-th iteration of the simplified ADMM algorithm for~\eqref{E:ADMMForm} proposed in~\cite{ODonoghue2016} consists of the following three steps, where $\mathbb{P}_{\mathcal{K}}$ denotes projection onto the cone $\mathcal K$:
\begin{subequations} \label{E:ADMMSteps}
    \begin{align}
    \label{E:HsdeStep1}
      \hat{u}^{(n+1)} &= (I+Q)^{-1}\left(u^{(n)}+ v^{(n)}\right), \\
    \label{E:HsdeStep2}
      u^{(n+1)} &= \mathbb{P}_{\mathcal{K}}\left(\hat{u}^{(n+1)}-v^{(n)}\right),\\
    \label{E:HsdeStep3}
      v^{(n+1)} &= v^{(n)} - \hat{u}^{(n+1)} + u^{(n+1)}.
    \end{align}
\end{subequations}

Note that~\eqref{E:HsdeStep2} is inexpensive, since $\mathcal{K}$ is the cartesian product of simple cones (zero, free and non-negative cones) and small PSD cones, and can be efficiently carried out in parallel. The third step is also computationally inexpensive and parallelizable. On the contrary, even when the preferred factorization of $I+Q$ (or its inverse) is cached before starting the iterations, a direct implementation of~\eqref{E:HsdeStep1} may require substantial computational effort because
$$
    Q \in \mathbb{S}^{n^2 + 2n_d + m + 1}
$$
is a very large matrix (\emph{e.g.}, $n^2 + 2n_d + m + 1 = 2\,360\,900$ for problem  \probname{rs365} in Section~\ref{se:simulationNonchordal}). Yet, it is evident from~\eqref{E:HSDEnotations} that $Q$ is highly structured and sparse, and these properties can be exploited to speed up step~\eqref{E:HsdeStep1} using a series of block-eliminations and the matrix inversion lemma~\cite[Section C.4.3]{boyd2004convex}.

\subsubsection{Solving the ``outer'' linear system}

The affine projection step~\eqref{E:HsdeStep1} requires the solution of a linear system (which we refer to as the ``outer'' system for reasons that will become clear below) of the form
\begin{equation} \label{E:LinearSystem}
\begin{bmatrix}
M & \zeta \\ -\zeta^T & 1
\end{bmatrix}
\begin{bmatrix}
\hat{u}_1 \\ \hat{u}_2
\end{bmatrix}
=
\begin{bmatrix}
\omega_1 \\ \omega_2
\end{bmatrix},
\end{equation}
where
\begin{align}
\label{E:MDef}
M &:= \begin{bmatrix}
I & -\hat{A}^T \\
\hat{A} & I
\end{bmatrix},
&
\zeta &:= \begin{bmatrix}\hat{c} \\ -\hat{b} \end{bmatrix},&
\hat{A} &:= \begin{bmatrix} A & 0 \\ H & -I \end{bmatrix}, &
\hat{c} &:=  \begin{bmatrix}c\\0 \end{bmatrix},&
\hat{b} &:= \begin{bmatrix}b\\0 \end{bmatrix}
\end{align}
and we have split
\begin{equation}
\label{E:OmegaDef}
u^{(n)}+ v^{(n)} = \begin{bmatrix}
\omega_1 \\ \omega_2
\end{bmatrix}.
\end{equation}
Note that $\hat{u}_2$ and $\omega_2$ are scalars. {Eliminating $\hat{u}_2$ from the first block equation in~\eqref{E:LinearSystem} yields}
\begin{subequations}
\begin{align}
\label{E:MatrixInverse}
(M+\zeta\zeta^T) \hat{u}_1 &= \omega_1 - \omega_2 \zeta,\\
\label{E:FirstBlockElimination}
\hat{u}_2 &= \omega_2 + \zeta^T \hat{u}_1.
\end{align}
\end{subequations}
Moreover, applying the matrix inversion lemma~\cite[Section C.4.3]{boyd2004convex} to~\eqref{E:MatrixInverse} shows that
\begin{equation} \label{E:MatrixInverseResult}
    \hat{u}_1  =
    \left[
    I -
    \frac{ (M^{-1} \zeta)\zeta^T}{1 + \zeta^T (M^{-1} \zeta)}
    \right]M^{-1}
    \left(
    \omega_1 - \omega_2 \zeta
    \right).
\end{equation}

Note that the vector $M^{-1} \zeta$ and the scalar $1 + \zeta^T (M^{-1} \zeta)$ depend only on the problem data, and can be computed before starting the ADMM iterations (since $M$ is quasi-definite it can be inverted, and any symmetric matrix obtained as a permutation of $M$ admits an LDL factorization). Instead, recalling from~\eqref{E:OmegaDef} that $\omega_1 - \omega_2 \zeta$ changes at each iteration because it depends on the iterates $u^{(n)}$ and $v^{(n)}$, the vector $M^{-1} \left(\omega_1 - \omega_2 \zeta\right)$ must be computed at each iteration. Consequently, computing $\hat{u}_1$ and $\hat{u}_2$ requires the solution of an ``inner'' linear system for the vector $M^{-1} \left(\omega_1 - \omega_2 \zeta\right)$, followed by inexpensive vector inner products and scalar-vector operations in~\eqref{E:MatrixInverseResult} and~\eqref{E:FirstBlockElimination}.

\subsubsection{Solving the ``inner'' linear system}

Recalling the definition of $M$ from~\eqref{E:MDef}, the ``inner'' linear system to calculate $\hat{u}_1$ in~\eqref{E:MatrixInverseResult} has the form
\begin{equation} \label{E:SecondBlockElimination}
    \begin{bmatrix}
    I & -\hat{A}^T \\
    \hat{A} & I
    \end{bmatrix}
    \begin{bmatrix}
    \sigma_1 \\ \sigma_2
    \end{bmatrix}
    = \begin{bmatrix}
    \nu_1 \\ \nu_2
    \end{bmatrix}.
\end{equation}
{Here}, $\sigma_1$ and $\sigma_2$ are the unknowns and represent suitable partitions of the vector $M^{-1}(\omega_1-\omega_2\zeta)$ in~\eqref{E:MatrixInverseResult}, which is to be calculated, and we have split
$$
\omega_1 - \omega_2 \zeta =
    \begin{bmatrix}
    \nu_1 \\ \nu_2
    \end{bmatrix}.
$$
Applying block elimination to remove $\sigma_1$ from the second equation in~\eqref{E:SecondBlockElimination}, we obtain
\begin{subequations}
\begin{align}
    (I+\hat{A}^T\hat{A})\sigma_1 &= \nu_1 + \hat{A}^T \nu_2, \label{E:SecBlkEliResult_s1}\\
    \sigma_2 &= -\hat{A} \sigma_1 + \nu_2. \label{E:SecBlkEliResult}
\end{align}
\end{subequations}
Recalling the definition of $\hat{A}$ and recognizing that
$$
D = H^T H= \sum_{k=1}^p H_k^TH_k
$$
is a diagonal matrix, as already noted in Section~\ref{se:MinYblkprimal}, we also have
\begin{equation*} 
I + \hat{A}^T\hat{A} = \begin{bmatrix} (I+D + A^TA) & -H^T \\ -H & 2I\end{bmatrix}.
\end{equation*}
Block elimination can therefore be used once again to solve \eqref{E:SecBlkEliResult_s1}, and simple algebraic manipulations show that the only matrix to be factorized (or inverted) is
\begin{equation} \label{E:FactorizationMatrix}
    I+\frac{1}{2}D + A^T A \in \mathbb{S}^{n^2}.
\end{equation}
Note that this matrix depends only on the problem data and the chordal decomposition, so it can be factorized/inverted before starting the ADMM iterations.
In addition, it is of the ``diagonal plus low rank" form because $A\in\mathbb{R}^{m\times n^2}$ with $m < n^2$ (in fact, often $m\ll n^2$). This means that the matrix inversion lemma can be used to reduce the size of the matrix to factorize/invert even further: letting $P = I+\frac{1}{2}D$ be the diagonal part of~\eqref{E:FactorizationMatrix}, we have
$$
    (P + A^T A)^{-1} = P^{-1} - P^{-1}A^T(I+AP^{-1}A^T)^{-1}AP^{-1}.
$$

In summary, after a series of block eliminations and applications of the matrix inversion lemma, step~\eqref{E:HsdeStep1} of the ADMM algorithm for~\eqref{E:ADMMForm} only requires the solution of an $m\times m$ linear system of equations with coefficient matrix
\begin{equation} \label{E:HSDEFactorizationFinal}
    I+A\left(I + \frac{1}{2}D\right)^{-1}A^T \in \mathbb{S}^{m},
\end{equation}
plus a sequence of matrix-vector, vector-vector, and scalar-vector multiplications. A detailed count of floating-point operations is given in Section~\ref{Section:Complexity}.

\subsubsection{Stopping conditions}
\label{s:hsde-stopping}

The ADMM algorithm described in the previous section can be stopped after the $n$-th iteration if a primal-dual optimal solution or a certificate of primal and/or dual infeasibility is found, up to a specified tolerance $\epsilon_\mathrm{tol}$.
As noted in~\cite{ODonoghue2016}, rather than checking the convergence of the variables $u$ and $v$, it is desirable to check the convergence of the original primal and dual SDP variables using the primal and dual residual error measures normally considered in interior-point algorithms~\cite{sturm1999using}. For this reason, we employ different  stopping conditions than those used in Algorithms~\ref{A:ADMMPrimal} and~\ref{A:ADMMDual}, which we define below {using the following notational convention: we denote the entries of $u$ and $v$ in~\eqref{E:HSDEnotations} that correspond to $x$, $y$, $\tau$, and $z$, respectively, by $u_x$, $u_y$, $u_{\tau}$, and $v_z$.}

If $u_{\tau}^{(n)} >0$ at the $n$-th iteration of the ADMM algorithm,  we take
\begin{align}
\label{e:hsde-solution}
    x^{(n)} = \frac{u_{x}^{(n)}}{u_{\tau}^{(n)}}, \qquad
    y^{(n)} = \frac{u_{y}^{(n)}}{u_{\tau}^{(n)}}, \qquad
    {z^{(n)} = \frac{H^Tv_{z}^{(n)}}{u_{\tau}^{(n)}}}
\end{align}
as the candidate primal-dual solutions, and define the relative primal residual, dual residual, and duality gap as
\begin{subequations}
	\label{E:HSDEerror}
	\begin{gather}
		\label{E:HSDEerror_primal}
    	\epsilon_\mathrm{p} :=  \frac{\|Ax^{(n)}-b\|_2}{1+\|b\|_2}, \\
    	\label{E:HSDEerror_dual}
    	\epsilon_\mathrm{d} :=  \frac{\|A^Ty^{(n)}+z^{(n)}-c\|_2}{1+\|c\|_2}, \\
    	\epsilon_\mathrm{g} :=  \frac{|c^Tx^{(n)}-b^Ty^{(n)}|}{1+\vert c^Tx^{(n)}\vert +\vert b^Ty^{(n)}\vert}.
	\end{gather}
\end{subequations}
Also, we define the residual in consensus constraints as
\begin{equation} \label{E:HSDEconsensusResidual}
    \epsilon_\mathrm{c} := \max \{\eqref{E:PrimalConsensusResidual},\eqref{E:DualConsensusResidual}\}.
\end{equation}
We terminate the algorithm if $\max \{\epsilon_\mathrm{p}, \epsilon_\mathrm{d}, \epsilon_\mathrm{g}, \epsilon_\mathrm{c}\}$ is smaller than $\epsilon_\mathrm{tol}$.  If $u_{\tau}^{(n)} =0$, instead, we terminate the algorithm if
\begin{equation} \label{E:Infeasibility}
 \max\left\{
        \|Au_x^{(n)}\|_2 + \frac{c^Tu_x^{(n)}}{\|c\|_2}\,\epsilon_\mathrm{tol},
        \,
         \|A^Tu_y^{(n)}+H^T v_z^{(n)}\|_2 - \frac{b^Tu_y^{(n)}}{\|b\|_2}\,\epsilon_\mathrm{tol}
\right\} \leq 0.
\end{equation}
Certificates of primal or dual infeasibility (with tolerance $\epsilon_\mathrm{tol}$) are then given, respectively, by the points
$u_y^{(n)} / ( b^T u_y^{(n)} )$ and $-u_x^{(n)}/ ( c^T u_x^{(n)} )$.
These stopping criteria are similar to those used by many other conic solvers, { and coincide with those used in SCS~\cite{scs} except for the addition of the residual in the consensus constraints~\eqref{E:HSDEconsensusResidual}.} The complete ADMM algorithm to solve the HSDE of the primal-dual pair of domain- and range-space decomposed SDPs is summarized in Algorithm~\ref{A:ADMMhsde}.

\begin{algorithm}[t]
	\caption{ADMM for the HSDE of sparse SDPs with chordal decomposition}
	\label{A:ADMMhsde}
	\begin{algorithmic}[1]
		\algsize
		\State Set $\epsilon_\mathrm{tol} >0$, a maximum  number of iterations $n_{\max}$ and initial guesses $\hat{u}^{(0)}$, $u^{(0)}$, $v^{(0)}$.
		\State Data preprocessing: chordal extension, chordal decomposition and factorization of the matrix in~\eqref{E:HSDEFactorizationFinal}.
		\For{$n = 1, \ldots, n_{\max}$}
		\State Compute $\hat{u}^{(n+1)}$ using the sequence of block eliminations~\eqref{E:LinearSystem}-\eqref{E:HSDEFactorizationFinal}.
		\State Compute $u^{(n+1)}$ using~\eqref{E:HsdeStep2}.
		\State Compute $v^{(n+1)}$ using~\eqref{E:HsdeStep3}.
		\If{$u_{\tau}^{(n)} >0$}
		\State Compute $\epsilon_\mathrm{p},\epsilon_\mathrm{d},{\epsilon_\mathrm{g}},{\epsilon_\mathrm{c}}$.
		\If{$\max\{\epsilon_\mathrm{p},\epsilon_\mathrm{d},{\epsilon_\mathrm{g}},{\epsilon_\mathrm{c}}\} \leq \epsilon_\mathrm{tol}$}
		\State \textbf{break}
		\EndIf
		\Else
		\If{ \eqref{E:Infeasibility} holds}
		\State \textbf{break}
		\EndIf
		\EndIf
		\EndFor
	\end{algorithmic}
\end{algorithm}

\section{Complexity analysis via flop count} \label{Section:Complexity}

The computational complexity of each iteration of Algorithms~\ref{A:ADMMPrimal}-\ref{A:ADMMhsde} can be assessed by counting the total number of {required} \emph{floating-point operations} {(flops)}---that is, the number of additions, subtractions, multiplications, or divisions of two floating-point numbers~\cite[Appendix C.1.1]{boyd2004convex}---as a function of problem dimensions. For~\eqref{E:OptCondMinYPrimal} and~\eqref{E:OptCondMinYDual} we have
$$
A \in \mathbb{R}^{m \times n^2}, \quad
b \in \mathbb{R}^{m}, \quad
c \in \mathbb{R}^{n^2}, \quad
D \in \mathbb{S}^{n^2},
\quad
H_k \in \mathbb{R}^{ |\mathcal{C}_k|^2 \times n^2} \text{ for } k = 1, \ldots, p,
$$
while the dimensions of the variables are 
$$
x \in \mathbb{R}^{n^2}, \quad
y \in \mathbb{R}^m, \quad
x_k,\,\lambda_k \in \mathbb{R}^{|\mathcal{C}_k|^2} \text{ for }  k = 1, \ldots, p.
$$

In this section, we count the flops in Algorithms~\ref{A:ADMMPrimal}--\ref{A:ADMMhsde} as a function of $m$, $n$, $p$, and $n_d = \sum_{k=1}^p|\mathcal{C}_k|^2$. We do not consider the sparsity in the problem data, both for simplicity and because sparsity is problem-dependent. Thus, the matrix-vector product $Ax$ is assumed to cost $(2n^2-1)m$ flops (for each row, we need $n^2$ multiplications and $n^2-1$ additions), while $A^Ty$ is assumed to cost $(2m-1)n^2$ flops. In practice, of course, these matrix-vector products may require significantly fewer flops if $A$ is sparse, and sparsity \textit{should} be exploited in any implementation to reduce computational cost. The only exception that we make concerns the matrix-vector products $H_k x$ and $H_k^Tx_k$ because each $H_k$, $k=1,\,\ldots,\,p$, is an ``entry-selector'' matrix that extracts the subvector $x_k \in \mathbb{R}^{|\mathcal{C}_k|^2}$ from $x\in\mathbb{R}^{n^2}$. Hence, the operations $H_k x$ and $H_k^Tx_k$ require no actual matrix multiplications but only indexing operations (plus, possibly, making copies of floating-point numbers depending on the implementation), so they cost no flops according to our definition. However, we do not take into account that the vectors $H_k^Tx_k \in \mathbb{R}^{n^2}$, $k=1,\,\ldots,\,p$, are often sparse, because their sparsity depends on the particular problem at hand. It follows from these considerations that computing the summation $\sum_{k=1}^pH_k^Tx_k$ costs $(p-1)n^2$ flops.

{
Using these rules, in the Appendix we prove the following results.
\begin{proposition} \label{T:FlopsLinP&D}
    Given the Cholesky factorization of $AD^{-1}A^T = LL^T$, where $L$ is lower triangular, solving the linear systems~\eqref{E:OptCondMinYPrimal} and~\eqref{E:OptCondMinYDual} via block elimination costs $(4m+p+3)n^2+2m^2+2n_d$  flops.
\end{proposition}
\begin{proposition}\label{T:FlopsLinHsde_final}
    Given the constant vector $\hat{\zeta} :=  (M^{-1} \zeta)/( 1 + \zeta^T M^{-1} \zeta ) \in \mathbb{R}^{n^2+2n_d+m}$ and the Cholesky factorization $I+A(I + \frac{1}{2}D)^{-1}A^T = LL^T$, where $L$ is lower triangular, solving~\eqref{E:HsdeStep1} using the sequence of block eliminations~\eqref{E:LinearSystem}--\eqref{E:HSDEFactorizationFinal} requires $(8m+2p+11)n^2+2m^2 + 7m +21n_d-1$ flops.
\end{proposition}
These propositions reveal that the computational complexity of the affine projections in Algorithms~\ref{A:ADMMPrimal} and~\ref{A:ADMMDual}, which amount to solving the linear systems~\eqref{E:OptCondMinYPrimal} and~\eqref{E:OptCondMinYDual}, is comparable to that of the affine projection~\eqref{E:HsdeStep1} in Algorithm~\ref{A:ADMMhsde}. In fact, since typically $m \ll n^2$, we expect that the affine projection step of Algorithm~\ref{A:ADMMhsde} will be only approximately  twice as expensive as the corresponding step in Algorithms~\ref{A:ADMMPrimal} and~\ref{A:ADMMDual} in terms of the number of flops, and therefore also in terms of CPU time (the numerical results presented in Table~\ref{T:ResultsLinProj}, Section~\ref{Section:BlockArrow} will confirm this expectation).}

Similarly, the following result (also proved in the Appendix) guarantees that the leading-order costs of the conic projections in Algorithms~\ref{A:ADMMPrimal}--\ref{A:ADMMhsde} are identical and, importantly, depend only on the size and number of the maximal cliques in the chordal decomposition, \textit{not} on the dimension $n$ of the original PSD cone in~\eqref{E:PrimalSDP}--\eqref{E:DualSDP}.
\begin{proposition} \label{T:FlopsConic}
	The computational costs of the conic projections in Algorithms~\ref{A:ADMMPrimal}--\ref{A:ADMMhsde} require $\mathcal{O}(\sum_{k=1}^p|\mathcal{C}_k|^3)$ floating-point operations.
\end{proposition}
In particular, the computational burden grows as a linear function of the number of cliques when their size is fixed, and as a {cubic} function of the clique size.

{Finally, we emphasize that Propositions~\ref{T:FlopsLinP&D}--\ref{T:FlopsConic} suggest that Algorithms~\ref{A:ADMMPrimal}--\ref{A:ADMMhsde} should solve a primal-dual pair of sparse SDPs more efficiently than the general-purpose ADMM method for conic programs of~\cite{ODonoghue2016}, irrespective of whether this is used before or after chordal decomposition. In the former case, the benefit comes from working with smaller PSD cones: one block-elimination in equation (28) of~\cite{ODonoghue2016} allows solving affine projection step~\eqref{E:HsdeStep1} in $\mathcal{O}(mn^2)$ flops, which is typically comparable to the flop count of Propositions~\ref{T:FlopsLinP&D} and~\ref{T:FlopsLinHsde_final},\footnote{{This can be seem more clearly after using the crude bound $n_d \leq pn^2$ in Propositions~\ref{T:FlopsLinP&D} and~\ref{T:FlopsLinHsde_final} and recalling that, for typical problems, $m \ll n^2$ and $p \ll n$.}} but the conic projection step costs $\mathcal{O}(n^3)$ flops, which for typical sparse SDPs is significantly larger than $\mathcal{O}(\sum_{k=1}^p|\mathcal{C}_k|^3)$. In the latter case, instead, the conic projection~\eqref{E:HsdeStep2} costs the same for all methods, but projecting the iterates onto the affine constraints becomes much more expensive according to our flop count when the sequences of block eliminations described in Section~\ref{sec:hsde} is not exploited fully.
}

\section{Implementation and numerical experiments} \label{sec:simulation}

We implemented Algorithms~\ref{A:ADMMPrimal}--\ref{A:ADMMhsde} in an open-source MATLAB solver which we call {{\CDCS}} (Cone Decomposition Conic Solver). We refer to our implementation of Algorithms~\ref{A:ADMMPrimal}--\ref{A:ADMMhsde} as {{\CDCS}}-primal, {{\CDCS}}-dual and {{\CDCS}}-hsde, respectively. This section briefly describes {{\CDCS}} and presents numerical results {on} sparse SDPs from SDPLIB~\cite{borchers1999sdplib}, {large} and sparse SDPs with nonchordal sparsity patterns from~\cite{andersen2010implementation}, and randomly generated SDPs with block-arrow sparsity pattern. Such problems have also been used as benchmarks in~\cite{andersen2010implementation,sun2014decomposition}.

In order to highlight the advantages of chordal decomposition, first-order algorithms, and their combination, the three algorithms in {{\CDCS}} are compared to the interior-point solver SeDuMi~\cite{sturm1999using}, and to the single-threaded direct implementation of the first-order algorithm of~\cite{ODonoghue2016} provided by the conic solver SCS~\cite{scs}. {All experiments were carried out on a PC with a 2.8 GHz Intel Core i7 CPU and 8GB of RAM and} the solvers were called  with termination tolerance $\epsilon_{\rm tol}=10^{-3}$, number of iterations limited to $2\,000$, and their default remaining parameters. The purpose of comparing {{\CDCS}} to a low-accuracy IPM is to demonstrate the advantages of combining FOMs with chordal decomposition, while a comparison to the high-performance first-order conic solver SCS highlights the advantages of chordal decomposition alone. {When possible,} accurate solutions ($\epsilon_{\rm tol}=10^{-8}$) were also computed using SeDuMi; these can be considered ``exact'', and used to assess how far the solution returned by {{\CDCS}} is from optimality. Note that tighter tolerances could be used also with {\CDCS} and SCS to obtain a more accurate solution, at the expense of increasing the number of iterations required to meet the convergence requirements. More precisely, given the proven convergence rate of general ADMM algorithms, any termination tolerance $\epsilon_{\rm tol}$ is {generally} reached in {at most} $\mathcal{O}(\frac{1}{\epsilon_{\rm tol}})$ iterations.

\subsection{{{\CDCS}}}

To the best of our knowledge, {{\CDCS}} is the first open-source first-order conic solver that exploits chordal decomposition for the PSD cones and is able to handle infeasible problems. Cartesian products of the following cones are supported: the cone of free variables $\mathbb R^n$, the non-negative orthant $\mathbb R^n_{+}$, second-order cones, and PSD cones. The current implementation is written in MATLAB and can be downloaded from\\[1ex]
\centerline{{ \small \url{https://github.com/oxfordcontrol/cdcs}}.}\\[1ex]
Note that although many steps of Algorithms~\ref{A:ADMMPrimal}--\ref{A:ADMMhsde} can be carried out in parallel, our implementation is sequential. Interfaces with the optimization toolboxes YALMIP ~\cite{lofberg2005yalmip} and SOSTOOLS~\cite{papachristodoulou2013sostools} are also available.

\subsubsection{Implementation details}

{{\CDCS}} applies chordal decomposition to all PSD cones. Following~\cite{vandenberghe2014chordal}, the sparsity pattern of each PSD cone is chordal extended using the MATLAB function {{\tt \small chol}} to compute a symbolic Cholesky factorization of the approximate minimum-degree permutation of the cone's adjacency matrix, returned by the MATLAB function {{\tt \small symamd}}. The maximal cliques of the chordal extension are then computed using a {{\tt \small.mex}} function from SparseCoLO~\cite{fujisawa2009user}.

As far as the steps of our ADMM algorithms are concerned, projections onto the PSD cone are performed using the MATLAB routine {\tt \small eig}, while projections onto other supported cones only use vector operations. The Cholesky factors of the $m \times m$ linear system coefficient matrix (permuted using {\tt \small symamd}) are cached before starting the ADMM iterations. The permuted linear system is solved at each iteration using the routines {\tt \small cs\_lsolve} and {\tt \small cs\_ltsolve} from the CSparse library~\cite{CSparse}.

{{\CDCS}} solves the decomposed problems~\eqref{E:DecomposedPrimalVector} and/or~\eqref{E:DecomposedDualVector} using any of Algorithms~\ref{A:ADMMPrimal}--\ref{A:ADMMhsde}, and then attempts to construct a primal-dual solution of the original SDPs~\eqref{E:PrimalSDP} and~\eqref{E:DualSDP} with a maximum determinant completion routine (see~\cite[Section 2]{fukuda2001exploiting},~\cite[Chapter 10.2]{vandenberghe2014chordal}) adapted from SparseCoLO~\cite{fujisawa2009user}.
{We adopted this approach for simplicity of implementation, even though we cannot guarantee that the principal sub-matrices $E_{\mathcal{C}_k} X^* E_{\mathcal{C}_k}^T$ of the partial matrix $X^*$ returned by {\CDCS} as the candidate solution are strictly positive definite (a requirement for the maximum determinant completion to exist). This may cause the current completion routine to fail, although for all cases in which we have observed failure, this was due to {\CDCS} returning a candidate solution with insufficient accuracy that was not actually PSD-completable. In any case, our current implementation issues a warning when the matrix completion routine fails; future versions of {\CDCS} will include alternative completion methods, such as that discussed in~\cite[Chapter 10.3]{vandenberghe2014chordal} and the minimum-rank PSD completion~\cite[Theorem 1]{Madani2015ADMM}, which work also in the lack of strict positive definiteness.
}

\subsubsection{Adaptive penalty strategy}

While the ADMM algorithms proposed in the previous sections converge independently of the choice of penalty parameter $\rho$, in practice its value strongly influences the number of iterations required for convergence. Unfortunately, analytic results for the optimal choice of $\rho$ are not available except for very special problems~\cite{ghadimi2015optimal,raghunathan2014alternating}. Consequently, in order to improve the convergence rate and make performance less dependent on the choice of $\rho$, {{\CDCS}} employs the dynamic adaptive rule. 
\begin{equation*}
\rho^{k+1} = \begin{cases}
\mu\,\rho^{(n)}
	& \text{if } \|\epsilon_\mathrm{p}^{(n)}\|_2
	\geq \nu \|\epsilon_\mathrm{d}^{(n)}\|_2,
\\
\mu^{-1}\rho^{(n)}
	& \text{if } \|\epsilon_\mathrm{d}^{(n)}\|_2
	\geq \nu \|\epsilon_\mathrm{p}^{(n)}\|_2,
\\
\rho^{(n)} &\text{otherwise.}
\end{cases}
\end{equation*}
Here, $\epsilon_\mathrm{p}^{(n)}$ and $\epsilon_\mathrm{d}^{(n)}$ are the primal and dual residuals at the $n$-th iteration, while $\mu$ and $\nu$ are parameters no smaller than 1. Note that since $\rho$ does not enter any of the matrices being factorized/inverted, updating its value is computationally {inexpensive}.

The idea of the rule above is to adapt $\rho$ to balance the convergence of the primal and dual residuals to zero; more details can be found in~\cite[Section 3.4.1]{boyd2011distributed}. Typical choices for the parameters (the default in {{\CDCS}}) are $\mu = 2$
and $\nu = 10$~\cite{boyd2011distributed}.

\subsubsection{Scaling the problem data}

The relative scaling of the problem data also affects the convergence rate of ADMM algorithms. {{\CDCS}} scales the problem data after the chordal decomposition step using a strategy similar to~\cite{ODonoghue2016}. In particular, the decomposed SDPs~\eqref{E:DecomposedPrimalVector} and~\eqref{E:DecomposedDualVector} can be rewritten as:
%
\refstepcounter{equation}\label{E:DecomposedPrimalScale}
\begin{align}
\tag{\theequation a,b}
  &\begin{aligned}
    \min_{\hat{x}} \quad & \hat{c}^T\hat{x} \\
    \text{subject to} \quad &  \hat{A}\hat{x} = \hat{b}    \\
                            & \hat{x} \in \mathbb{R}^{n^2} \times \mathcal{K},
  \end{aligned}
  &
  &\begin{aligned}
        \max_{\hat{y}, \hat{z}} \qquad &  \hat{b}^T \hat{y} \\
        \text{subject to} \quad & \hat{A}^T\hat{y} + \hat{z} = \hat{c} \\
        & \hat{z} \in \{0\}^{n^2} \times \hat{\mathcal{K}}^*
    \end{aligned}
\end{align}
where
\begin{align*}
  \hat{x} &= \left[\begin{array}{c} x \\ s \\\end{array}\right], &
  \hat{y} &= \left[\begin{array}{c} y \\ t \\\end{array}\right], &
  \hat{z} &= \left[\begin{array}{c} 0 \\ z \\\end{array}\right], &
  \hat{c} &= \left[\begin{array}{c} c \\ 0 \\\end{array}\right], &
  \hat{b} &= \left[\begin{array}{c} b \\ 0 \\\end{array}\right], &
  \hat{A} &= \begin{bmatrix} A & 0 \\H & -I \\\end{bmatrix} .
\end{align*}
{{\CDCS}} solves the scaled {decomposed} problems
%
\refstepcounter{equation}\label{E:ScaledSDP}
\begin{align}
\tag{\theequation a,b}
  &\begin{aligned}
        \min_{\hat{x}} \qquad & \sigma (D\hat{c})^T \bar{x} \\
        \text{subject to} \quad & E\hat{A}D\bar{x} = \rho E\hat{b} \\
        & \bar{x} \in \mathbb{R}^{n^2} \times \mathcal{K},
    \end{aligned}
    &
    &\begin{aligned}
        \max_{\hat{y}, \hat{z}} \qquad & \rho (Eb)^T \bar{y} \\
        \text{subject to} \quad & D\hat{A}^TE\bar{y} + \bar{z} = \sigma D \hat{c} \\
        & \bar{z} \in \{0\}^{n^2} \times \mathcal{K}^*,
    \end{aligned}
\end{align}
obtained by scaling vectors $\hat{b}$ and $\hat{c}$ by positive scalars $\rho$ and $\sigma$, and the primal and dual equality constraints by positive definite, diagonal matrices $D$ and $E$.  Note that such a rescaling does not change the sparsity pattern of the problem. As already observed in~\cite{ODonoghue2016}, a good choice for $E$, $D$, $\sigma$ and $\rho$ is such that the rows of $\bar{A}$ and $\bar{b}$ have Euclidean norm close to one, and the columns of $\bar{A}$ and $\bar{c}$ have similar norms. If $D$ and $D^{-1}$ are chosen to preserve membership to the cone $\mathbb{R}^{n^2}\times\mathcal{K}$ and its dual, respectively (how this can be done is explained in~\cite[Section 5]{ODonoghue2016}), an optimal point for~\eqref{E:DecomposedPrimalScale} can be recovered from the solution of~\eqref{E:ScaledSDP}:
$$
    \hat{x}^* = \frac{D\bar{x}^*}{\rho}, \quad \hat{y}^* = \frac{E\bar{y}^*}{\sigma}, \quad \hat{z}^* = \frac{D^{-1}\bar{z}^*}{\sigma}.
$$

\subsection{Sparse SDPs from SDPLIB}

Our first experiment is based on large-scale benchmark problems from SDPLIB~\cite{borchers1999sdplib}: two Lov\'asz $\vartheta$ number SDPs (\probname{theta1} and \probname{theta2}), two infeasible SDPs (\probname{infd1} and \probname{infd2}), two MaxCut problems (\probname{maxG11} and \probname{maxG32}), and two SDP relaxations of box-constrained quadratic programs (\probname{qpG11} and \probname{qpG51}). Table~\ref{T:SparseStatistic} reports the dimensions of these problems, as well as chordal decomposition details.  Problems \probname{theta1} and \probname{theta2} are dense, so have only one maximal clique; all other problems are sparse and have many maximal cliques of size much smaller than the original cone.

\begin{table*}
    \centering
    \renewcommand\arraystretch{0.9}
    \caption{Details of the SDPLIB problems considered in this work. }
    \label{T:SparseStatistic}
\begin{tabular*}{\linewidth}{@{\extracolsep{\fill}}r rr r rr r rrrr}
    \toprule[1pt]
    & \multicolumn{2}{c}{Small} &
    & \multicolumn{2}{c}{Infeasible} &
    & \multicolumn{4}{c}{Large and sparse} \\
    \cline{2-3} \cline{5-6} \cline{8-11}\\[-0.5em]
    & {\tt\scriptsize theta1} & {\tt\scriptsize theta2} &
     & {\tt\scriptsize infd1} & {\tt\scriptsize infd2}  &
     & {\tt\scriptsize maxG11} & {\tt\scriptsize maxG32} &
       {\tt\scriptsize qpG11} & {\tt\scriptsize qpG51}    \\
    \cline{2-3} \cline{5-6} \cline{8-11}\\[-0.5em]
Original cone size, $n$  & 50  & 100 && 30  & 30  && 800 & 2\,000 & 1\,600 & 2\,000   \\
Affine constraints, $m$  & 104 & 498 && 10 & 10 && 800 & 2\,000 & 800  & 1\,000  \\
Number of cliques, $p$  & 1   & 1   && 1   & 1   && 598 & 1\,499 & 1\,405 & 1\,675  \\
Maximum clique size   & 50  & 100 && 30  & 30  && 24  & 60   & 24   & 304   \\
Minimum clique size   & 50  & 100 && 30  & 30  && 5   & 5    & 1    & 1      \\
        \bottomrule[1pt]
        \end{tabular*}
\end{table*}

\begin{table*}
	\centering
	\renewcommand\arraystretch{0.9}
	\caption{Results for two small SDPs, \probname{theta1} and \probname{theta2}, in SDPLIB.}
	\label{T:ResultsSmall}
	\begin{tabular*}{\linewidth}{@{\extracolsep{\fill}}r rrr r rrr}
		\toprule[1pt]
		& \multicolumn{3}{c}{{\tt\scriptsize theta1}} && \multicolumn{3}{c}{{\tt\scriptsize theta2}} \\
		\cline{2-4} \cline{6-8}\\[-0.5em]
		& Time (s) & \# Iter. & Objective &
		& Time (s) & \# Iter. & Objective  \\
		\cline{2-4} \cline{6-8}\\[-0.5em]
		SeDuMi (high) & 0.281 & 14 & 23.00 &
		& 1.216 & 15 & 32.88  \\
		SeDuMi (low) & 0.161 & 8 & 23.00 &
		& 0.650  & 8 & 32.88 \\
		SCS (direct) & 0.057 & 140 & 22.99&
		& 0.244 & 200 & 32.89 \\
		{{\CDCS}}-primal & 0.297 & 163 & 22.92  &
		& 0.618 & 188  & 32.94\\
		{{\CDCS}}-dual & 0.284 & 154  & 22.83 &
		& 0,605 &  178 &32.89  \\
		{{\CDCS}}-hsde & 0.230 & 156 & 23.03 &
		&  0.392 & 118  & 32.88 \\
		\bottomrule[1pt]
	\end{tabular*}
\end{table*}

The numerical results are summarized in Tables~\ref{T:ResultsSmall}--\ref{T:ResultsCPUIteration}. Table~\ref{T:ResultsSmall} shows that the small dense SDPs \probname{theta1} and \probname{theta2}, were solved in approximately the same CPU time by all solvers. Note that since these problems only have one maximal clique, SCS and {{\CDCS}}-hsde use similar algorithms, and performance differences are mainly due to the implementation (most notably, SCS is written in C). Table~\ref{T:ResultsInf} confirms that {{\CDCS}}-hsde successfully detects infeasible problems, while {{\CDCS}}-primal and {{\CDCS}}-dual do not have this ability.

The CPU time, number of iterations and terminal objective value for the four large-scale sparse SDPs {\probname{maxG11}, \probname{maxG32}, \probname{qpG11} and \probname{qpG51}} are listed in Table~\ref{T:ResultsMaxCut}. All algorithms in {{\CDCS}} were faster than either SeDuMi or SCS, especially for problems with smaller maximum clique size as one would expect {in light of the complexity analysis of Section~\ref{Section:Complexity}}. Notably, {{\CDCS}} solved \probname{maxG11}, \probname{maxG32}, and \probname{qpG11} in less than $100\,\rm s$, a speedup of approximately $9 \times$, $43 \times$, and $66\times$ over SCS. In addition, even though FOMs are only meant to provide moderately accurate solutions, the terminal objective value returned by {{\CDCS}}-hsde was always within 0.2\% of the high-accuracy optimal value computed using SeDuMi. This is an acceptable difference in many practical applications.

%
\begin{table*}
	\centering
	\renewcommand\arraystretch{0.9}
	\caption{Results for two infeasible SDPs in SDPLIB. An objective value of +Inf denotes infeasiblity. Results for the primal-only and dual-only algorithms in {{\CDCS}} are not reported since they cannot detect infeasibility.}
	\label{T:ResultsInf}
	\begin{tabular*}{\linewidth}{@{\extracolsep{\fill}}r rrr c rrr}
		\toprule[1pt]
		& \multicolumn{3}{c}{{\tt\scriptsize infp1}} && \multicolumn{3}{c}{{\tt\scriptsize infp2}} \\
		\cline{2-4} \cline{6-8}\\[-0.5em]
		& Time (s) & \# Iter. & Objective &
		& Time (s) & \# Iter. & Objective  \\
		\cline{2-4} \cline{6-8}\\[-0.5em]
		SeDuMi (high) & 0.127 &  2 &   +Inf && 0.033 & 2 & +Inf\\
		SeDuMi (low) &  0.120  &  2 &  +Inf &&  0.031 &  2 & +Inf\\
		SCS (direct) & 0.067 & 20 &   +Inf &&  0.031 & 20 & +Inf\\
		{{\CDCS}}-hsde &  0.109 &  118 & +Inf && 0.114&  101 & +Inf\\
		\bottomrule[1pt]
	\end{tabular*}
\end{table*}

\begin{table*}
	\centering
	\renewcommand\arraystretch{0.9}
	\caption{Results for four large-scale sparse SDPs in SDPLIB, \probname{maxG11}, \probname{maxG32}, \probname{qpG11} and \probname{qpG51}.}
	\label{T:ResultsMaxCut}
	\begin{tabular*}{\linewidth}{@{\extracolsep{\fill}}r rrr c rrr}
		\toprule[1pt]
		& \multicolumn{3}{c}{{\tt\scriptsize maxG11}} && \multicolumn{3}{c}{{\tt\scriptsize maxG32}} \\
		\cline{2-4} \cline{6-8}\\[-0.5em]
		& Time (s) & \# Iter. & Objective &
		& Time (s) & \# Iter. & Objective  \\
		\cline{2-4} \cline{6-8}\\[-0.5em]
		SeDuMi (high) &  88.9 & 13 & 629.2 &  &1\,266   & 14 &  1\,568\\
		SeDuMi (low) &  48.7 &  7 & 628.7 &  &  624  &  7  & 1\,566\\
		SCS (direct) &  93.9 &  1\,080 &  629.1 & &  2\,433   & 2\,000  & 1\,568 \\
		{{\CDCS}}-primal & 22.2   &  230  &  629.5  &  	& 84    &  311  & 1\,569  \\
		{{\CDCS}}-dual &  16.9  & 220 & 629.2 &  &  61  &  205  &1\,567\\
		{{\CDCS}}-hsde &  10.9 &  182 &  629.3 & &  56  &  291   & 1\,568\\
		\midrule
		& \multicolumn{3}{c}{{\tt\scriptsize qpG11}} && \multicolumn{3}{c}{{\tt\scriptsize qpG51}} \\
		\cline{2-4} \cline{6-8}\\[-0.5em]
		& Time (s) & \# Iter. & Objective &
		& Time (s) & \# Iter. & Objective  \\
		\cline{2-4} \cline{6-8}\\[-0.5em]
		SeDuMi (high) &   650 &  14 &  2\,449 &
		&  1\,895 & 22 & 1\,182 \\
		SeDuMi (low) &  357  &  8 &  2\,448 &
		& 1\,530  & 18  &  1\,182  \\
		SCS (direct) &   1\,065    &  2\,000 &  2\,449 &
		&  2\,220  &  2\,000 & 1\,288 \\
		{{\CDCS}}-primal &  29  &  249 &  2\,450  &
		&   482  &  1\,079 & 1\,145  \\
		{{\CDCS}}-dual & 21 &  193  &  2\,448 &
		&  396 &  797  & 1\,201  \\
		{{\CDCS}}-hsde &  16 &  219  &  2\,449 &
		&  865  & 2\,000 & 1\,182  \\
		\bottomrule[1pt]
	\end{tabular*}
\end{table*}

{
To provide further evidence to assess the relative performance of the tested solvers, Tables~\ref{T:ResultsErrors} and~\ref{T:SedumiSCSresiduals} report the constraint violations for the original (not decomposed) SDPs, alonside the error in the consensus constraints for the decomposed problems. Specifically:
\begin{enumerate}
	\item For {\CDCS}-primal, we measure how far the partial matrix $X = \mat(x) \in \mathbb{S}^n(\mathcal{E},?)$ is from being PSD-completable. This is the only quantity of interest because the equality constraints in~\eqref{E:PrimalSDP} are satisfied exactly by virtue of the second block equation in~\eqref{E:OptCondMinYPrimal}. Instead of calculating the distance between $X$ and the cone $\mathbb{S}^n_+(\mathcal{E},?)$ exactly using, for instance, the methods of~\cite{sun2015decomposition}) we bound it from above by computing the smallest non-negative constant ${\alpha}$ such that $X+{\alpha} I \in \mathbb{S}^n_+(\mathcal{E},?)$; {indeed, for such $\alpha$ it is clear that $\min_{Y\in\mathbb{S}^n_+(\mathcal{E},0)} \|Y - X\|_F \leq \| (X+{\alpha} I) - X\|_F = \alpha \sqrt{n}$.}
	This strategy is more economical because, letting $\lambda_{\rm min}(M)$ be the minimum eigenvalue of a matrix $M$, Theorem~\ref{T:ChordalCompletionTheorem} implies that
	$$\alpha = -\min \left\{0,\, \lambda_{\rm min} \left[ \mat\left(H_1 x \right) \right],\,\ldots,\,\lambda_{\rm min} \left[ \mat\left(H_p x \right) \right] \right\}.$$
	To mitigate the dependence on the scaling of $X$, Table~\ref{T:ResultsErrors} lists the normalized error
	\begin{equation}\label{e:PSD-violation}
	\epsilon_{\alpha} := \frac{\alpha}{1+\|X\|_F}.
	\end{equation}	

	\item For {\CDCS}-dual, given the candidate solutions $y$ and $Z = \mat(\sum_{k=1}^p H_k^T z_k) =  \sum_{k=1}^p E_{\mathcal{C}_k} \mat(z_k) E_{\mathcal{C}_k}^T$ we report the violation of the equality constraints in~\eqref{E:DualSDP} given by the relative dual residual $\epsilon_{\mathrm{d}}$, defined as in~\eqref{E:HSDEerror_dual}. {Note that~\eqref{E:zkUpdate}  guarantees that the matrices $\mat(z_1),\,\ldots,\,\mat(z_p)$ are PSD, so $Z$ is also PSD.}
	
	\item For the candidate solution~\eqref{e:hsde-solution} returned by {\CDCS}-hsde, we list the relative primal and dual residuals $\epsilon_{\rm p}$ and $\epsilon_{\rm d}$ defined in~\eqref{E:HSDEerror_primal} and~\eqref{E:HSDEerror_dual}, as well as the error measure $\epsilon_\alpha$ computed with~\eqref{e:PSD-violation}.
	
	\item For SeDuMi and SCS, we report only the primal and dual residuals~\eqref{E:HSDEerror_primal}--\eqref{E:HSDEerror_dual} since the PSD constraints are automatically satisfied in both the primal and the dual problems.
\end{enumerate}

\begin{table*}[t]
	\centering
	\setlength{\abovecaptionskip}{1pt}
	\setlength{\belowcaptionskip}{0em}
	\renewcommand\arraystretch{0.9}
	\caption{{Residuals for the solutions returned by {\CDCS} with $\epsilon_{\text{tol}} = 10^{-3}$
			and the maximum number of iterations fixed to $2\,000$. The residual $\epsilon_{\mathrm{c}}$, defined in~\eqref{E:PrimalConsensusResidual}, \eqref{E:DualConsensusResidual}, and~\eqref{E:HSDEconsensusResidual} for {\CDCS}-primal, {\CDCS}-dual, and {\CDCS}-hsde respectively, measures the error in the consensus constraints
			of the decomposed SDPs. The quantities $\epsilon_{\mathrm{p}}$ and $\epsilon_{\mathrm{d}}$ are defined in~\eqref{E:HSDEerror_primal} and~\eqref{E:HSDEerror_dual}, respectively, and measure the primal-dual residuals of the equality constraints for the original SDPs (before decomposition). Finally, $\epsilon_\alpha$ computed using~\eqref{e:PSD-violation} measures the violation of the semidefiniteness constraint in the primal SDP.}
	}
	\label{T:ResultsErrors}
	\begin{tabular*}{\linewidth}{@{\extracolsep{\fill}}l c cc c cc c cccc}
		\toprule[1pt]
		&  \multicolumn{2}{c}{{\CDCS}-primal} && \multicolumn{2}{c}{{\CDCS}-dual} &&  \multicolumn{4}{c}{{\CDCS}-hsde}  \\
		\cline{2-3} \cline{5-6} \cline{8-11} \\[-0.8em]
		& $\epsilon_{\mathrm{c}}$ & $\epsilon_\alpha$
		&& $\epsilon_{\mathrm{c}}$ & $\epsilon_{\mathrm{d}}$
		&& $\epsilon_\alpha$ & $\epsilon_{\mathrm{p}}$ &$\epsilon_{\mathrm{d}}$ & $\epsilon_{\mathrm{c}}$ \\
		\cline{2-3} \cline{5-6} \cline{8-11} \\[-0.8em]
		{\tt\small maxG11} & $9.97$e-$4$ & $2.72$e-$4$ && $1.34$e-$4$ & $2.48$e-$4$ &&  $2.84$e-$4$ & $5.52$e-$4$ & $3.46$e-$4$ & $9.98$e-$4$ \\
		{\tt\small maxG32} & $9.96$e-$4$ & $2.22$e-$4$ && $4.75$e-$4$ & $11.1$e-$4$ &&  $1.78$e-$4$ & $2.74$e-$4$ & $2.68$e-$4$ & $9.97$e-$4$ \\
		{\tt\small qpG11}  & $3.70$e-$4$ & $3.62$e-$4$ && $9.94$e-$4$ & $20.0$e-$4$ &&  $2.76$e-$4$ & $9.66$e-$4$ & $3.60$e-$4$ & $8.40$e-$4$ \\
		{\tt\small qpG51}  & $9.34$e-$4$ & $0.54$e-$4$ && $9.91$e-$4$ & $4.02$e-$4$ &&  $0.15$e-$4$ & $13.3$e-$4$ & $31.3$e-$4$ & $2.10$e-$4$ \\
		{\tt\small rs35}   & $9.99$e-$4$ & $7.05$e-$4$ && $6.41$e-$4$ & $9.36$e-$4$ &&  $0.01$e-$4$ & $0.31$e-$4$ & $13.6$e-$4$ & $9.99$e-$4$ \\
		{\tt\small rs200}  & $9.92$e-$4$ & $6.73$e-$4$ && $1.48$e-$4$ & $2.47$e-$4$ &&  $2.74$e-$4$ & $9.37$e-$4$ & $9.95$e-$4$ & $9.71$e-$4$ \\
		{\tt\small rs228}  & $9.97$e-$4$ & $5.92$e-$4$ && $2.14$e-$4$ & $2.71$e-$4$ &&  $1.75$e-$4$ & $8.82$e-$4$ & $9.83$e-$4$ & $9.86$e-$4$ \\
		{\tt\small rs365}  & $9.98$e-$4$ & $2.98$e-$4$ && $4.61$e-$4$ & $9.99$e-$4$ &&  $0.36$e-$4$ & $2.50$e-$4$ & $9.99$e-$4$ & $9.97$e-$4$ \\
		{\tt\small rs1555} & $2.21$e-$4$ & $1.55$e-$4$ && $9.97$e-$4$ & $19.6$e-$4$ &&  $0.98$e-$4$ & $11.0$e-$4$ & $1.89$e-$4$ & $9.99$e-$4$ \\
		{\tt\small rs1907} & $9.98$e-$4$ & $3.27$e-$4$ && $6.24$e-$4$ & $12.0$e-$4$ &&  $0.23$e-$4$ & $3.78$e-$4$ & $9.96$e-$4$ & $9.90$e-$4$ \\
		\bottomrule[1pt]
	\end{tabular*}
\end{table*}

\begin{table*}[t]
	\centering
	\setlength{\abovecaptionskip}{1pt}
	\setlength{\belowcaptionskip}{0em}
	\renewcommand\arraystretch{0.9}
	\caption{{Residuals $\epsilon_{\mathrm{p}}$ and $\epsilon_{\mathrm{d}}$, defined as in~\eqref{E:HSDEerror_primal} and~\eqref{E:HSDEerror_dual}, for the solutions returned by SeDuMi and SCS. Both solvers were called with $\epsilon_{\text{tol}} = 10^{-3}$, and the maximum number of iterations for SCS is $2\,000$. Entries marked *** denote failure due to memory limitations.}
	}
	\label{T:SedumiSCSresiduals}
	\begin{tabular*}{0.54\linewidth}{l c cc c c}
		\toprule[1pt]
		&  \multicolumn{2}{c}{SeDuMi} && \multicolumn{2}{c}{SCS} \\
		\cline{2-3} \cline{5-6} \\[-0.8em]
		& $\epsilon_{\mathrm{p}}$ & $\epsilon_{\mathrm{d}}$
		&& $\epsilon_{\mathrm{p}}$ & $\epsilon_{\mathrm{d}}$ \\
		\cline{2-3} \cline{5-6} \\[-0.8em]
		{\tt\small maxG11} & $8.36$e-$6$ & $5.95$e-$7$ && $4.71$e-$6$ & $9.98$e-$4$ \\
		{\tt\small maxG32} & $8.45$e-$6$ & $3.85$e-$7$ && $2.46$e-$6$ & $3.20$e-$3$ \\
		{\tt\small qpG11}  & $1.13$e-$5$ & $2.90$e-$7$ && $1.32$e-$5$ & $2.85$e-$2$ \\
		{\tt\small qpG51}  & $4.23$e-$6$ & $5.41$e-$7$ && $1.15$e-$4$ & $2.36$e-$1$ \\
		{\tt\small rs35}   & $4.33$e-$7$ & $1.73$e-$10$ && $1.41$e-$5$ & $7.60$e-$2$ \\
		{\tt\small rs200}  & $5.41$e-$6$ & $4.83$e-$9$ && $1.55$e-$5$ & $4.43$e-$1$ \\
		{\tt\small rs228}  & $3.40$e-$6$ & $1.89$e-$9$ && $2.09$e-$5$ & $1.72$e-$1$ \\
		{\tt\small rs365}  & *** & *** && $2.70$e-$5$ & $7.20$e-$1$ \\
		{\tt\small rs1555} & *** & *** && $3.67$e-$6$ & $9.22$e-$1$ \\
		{\tt\small rs1907} & *** & *** && $3.29$e-$5$ & $7.96$e-$1$ \\
		\bottomrule[1pt]
	\end{tabular*}
\end{table*}

The results in Table~\ref{T:ResultsErrors} demonstrate that, for the problems tested in this work, the residuals for the original SDPs are comparable to the convergence tolerance used in {\CDCS} even when they are not tracked directly. The performance of SCS on our test problems is relatively poor. It is well known that the performance of ADMM algorithms is sensitive to their parameters, as well as problem scaling. We have calibrated CDSC using typical parameter values that offer a good compromise between efficiency and reliability, but we have not tried to fine-tune SCS under the assumption that good parameter values have already been chosen by its developers. Although performance may be improved through further parameter optimization, the discrepancy between the primal and dual residuals reported in Table~\ref{T:SedumiSCSresiduals} suggests that slow convergence may be due to problem scaling for these instances. {Note that {\CDCS} and SCS adopt the same rescaling strategy, but one key difference is that {\CDCS} applies it to the decomposed SDP rather than to the original one. Thus, CDCS has more degrees of scaling freedom than SCS, which might be the reason for the substantial improvement in convergence performance.} Further investigation of the effect of scaling in ADMM-based algorithms for conic programming, however, is beyond the scope of this work.
}

Finally, to offer a comparison of the performance of {\CDCS} and SCS that is insensitive both to problem scaling and to differences in the stopping conditions, Table~\ref{T:ResultsCPUIteration} reports the average CPU time per iteration required to solve the sparse SDPs \probname{maxG11}, \probname{maxG32}, \probname{qpG11} and \probname{qpG51}, as well as the dense SDPs \probname{theta1} and \probname{theta2}. Evidently, all algorithms in {{\CDCS}} are faster than SCS for the large-scale sparse SDPs (\probname{maxG11}, \probname{maxG32}, \probname{qpG11} and \probname{qpG51}), and in particular {{\CDCS}}-hsde improves on SCS by approximately $1.8 \times$, $8.7 \times$, $8.3 \times$, and $2.6 \times$ for each problem, respectively. This is to be expected since the conic projection step in {{\CDCS}} is more efficient due to smaller semidefinite cones, but the results are remarkable considering that {{\CDCS}} is written in MATLAB, while SCS is implemented in C. Additionally, the performance of {{\CDCS}} could be improved even further with a parallel implementation of the projections onto small PSD cones.

\begin{table}
\centering
\renewcommand\arraystretch{0.9}
\caption{Average CPU time per iteration (in seconds) for the SDPs from SDPLIB  tested in this work.}
\label{T:ResultsCPUIteration}
\begin{tabular*}{\linewidth}{@{\extracolsep{\fill}}r rrrrrr}
\toprule[1pt]
 & {\tt\scriptsize theta1} & {\tt\scriptsize theta2}
 & {\tt\scriptsize maxG11} & {\tt\scriptsize maxG32}
 & {\tt\scriptsize qpG11} & {\tt\scriptsize qpG51}\\
\cline{2-7}\\[-0.5em]
SCS (direct) & 4.0$\times 10^{-4}$ & 1.2$\times 10^{-3}$ & 0.087 & 1.216 & 0.532 & 1.110 \\
{{\CDCS}}-primal & 1.8$\times 10^{-3}$ & 3.3$\times 10^{-3}$ &  0.076  & 0.188  & 0.101 & 0.437  \\
{{\CDCS}}-dual & 1.8$\times 10^{-3}$ & 3.4$\times 10^{-3}$ & 0.064  & 0.174  & 0.091 & 0.484\\
{{\CDCS}}-hsde &   1.5$\times10^{-3}$ & 3.3$\times 10^{-3}$ &  0.048 & 0.140 & 0.064 & 0.430\\
\bottomrule[1pt]
\end{tabular*}
\end{table}

\subsection{Nonchordal SDPs} \label{se:simulationNonchordal}

\begin{figure}
    \centering
    \setlength{\abovecaptionskip}{0pt}
    \setlength{\belowcaptionskip}{0em}
    \newcommand{\fighspace}{\hspace{1cm}}
    \subfigure[\probname{rs35}]
    { \label{fig:nonchordal_a}
		\includegraphics[scale=0.9]{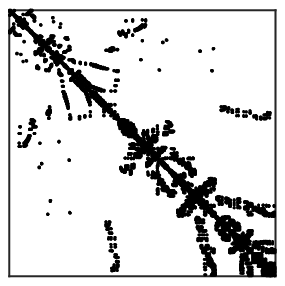}
    }
    \fighspace
    \subfigure[\probname{ rs200} ]
    { \label{fig:nonchordal_b}
		\includegraphics[scale=0.9]{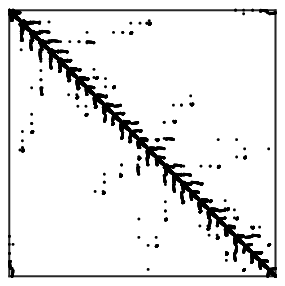}
    }
    \fighspace
    \subfigure[\probname{ rs228} ]
    { \label{fig:nonchordal_c}
		\includegraphics[scale=0.9]{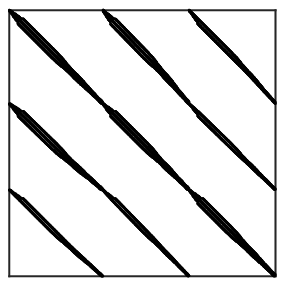}
    }
    \\[-1em]
    \subfigure[\probname{ rs365} ]
    { \label{fig:nonchordal_d}
		\includegraphics[scale=0.9]{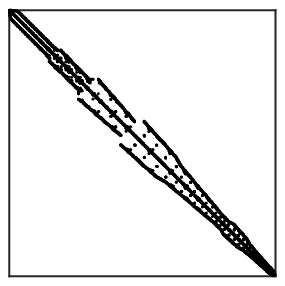}
    }
    \fighspace
    \subfigure[\probname{ rs1555} ]
    { \label{fig:nonchordal_e}
		\includegraphics[scale=0.9]{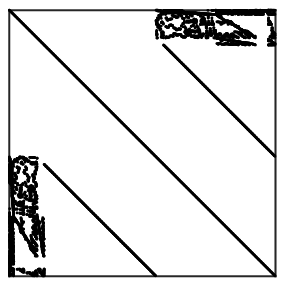}
    }
    \fighspace
    \subfigure[\probname{ rs1907} ]
    { \label{fig:nonchordal_f}
		\includegraphics[scale=0.9]{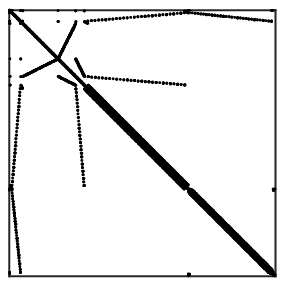}
    }
	\caption{Aggregate sparsity patterns of the nonchordal SDPs in~\cite{andersen2010implementation}; see Table~\ref{T:NonChordalStatistic} for the matrix dimensions.}
    \label{F:SP_Nonchordal}
\end{figure}

\begin{table*}[t]
    \centering
    \renewcommand\arraystretch{0.9}
    \caption{Summary of chordal decomposition for the chordal extensions  of the nonchordal SDPs form~\cite{andersen2010implementation}.
    }
    \label{T:NonChordalStatistic}
\begin{tabular*}{\linewidth}{@{\extracolsep{\fill}}r rrrrrr}
        \toprule[1pt]
        & {\tt\scriptsize rs35} & {\tt\scriptsize rs200}
         & {\tt\scriptsize rs228} & {\tt\scriptsize rs365}
         & {\tt\scriptsize rs1555} & {\tt\scriptsize rs1907} \\
        \cline{2-7}\\[-0.5em]
        Original cone size, $n$  & 2003& 3025 & 1919  & 4704  & 7479 & 5357   \\
        Affine constraints, $m$  & 200 & 200  & 200 & 200 & 200 & 200  \\
        Number of cliques, $p$   & 588 & 1635 & 783   & 1244   & 6912 & 611  \\
        Maximum clique size      & 418  & 102 & 92  & 322  & 187  & 285   \\
        Minimum clique size      & 5  & 4 & 3  & 6  & 2  & 7         \\
        \bottomrule[1pt]
        \end{tabular*}
\end{table*}

In our second experiment, we solved six large-scale SDPs with nonchordal sparsity patterns form~\cite{andersen2010implementation}: \probname{rs35}, \probname{rs200}, \probname{rs228}, \probname{rs365}, \probname{rs1555}, and \probname{rs1907}. The aggregate sparsity patterns of these problems, illustrated in Fig.~\ref{F:SP_Nonchordal}, come from the University of Florida Sparse Matrix Collection~\cite{davis2011university}. Table~\ref{T:NonChordalStatistic} demonstrates that all six sparsity patterns admit chordal extensions with maximum cliques that are much smaller than the original cone.

\begin{table*}
\centering
\renewcommand\arraystretch{0.9}
\caption{Results for large-scale SDPs with nonchordal sparsity patterns form~\cite{andersen2010implementation}. { Entries marked {***} indicate that the problem could not be solved due to memory limitations.}}
\label{T:ResultsNonChordal}
\begin{tabular*}{\linewidth}{@{\extracolsep{\fill}}r rrr c rrr}
\toprule[1pt]
& \multicolumn{3}{c}{{\tt\scriptsize rs35}} && \multicolumn{3}{c}{{\tt\scriptsize rs200}} \\
\cline{2-4} \cline{6-8}\\[-0.5em]
& Time (s) & \# Iter. & Objective &
& Time (s) & \# Iter. & Objective  \\
\cline{2-4} \cline{6-8}\\[-0.5em]
SeDuMi (high) &   1\,391 & 17 & 25.33 &
		&  4\,451 & 17 &  99.74   \\
SeDuMi (low)  & 986 & 11 & 25.34 &
					& 2\,223 & 8 & 99.73     \\
SCS (direct) & 2\,378  & 2\,000 &25.08 &
	& 9\,697 & 2\,000  & 81.87  \\
{{\CDCS}}-primal & 370 & 379 & 25.27 & & 159 & 577 & 99.61  \\
{{\CDCS}}-dual  & 272 & 245 &25.53 & & 103  & 353   &99.72 \\
{{\CDCS}}-hsde  & 2\,019 & 2\,000 & 25.47 & & 254  & 1\,114  & 99.70  \\
\midrule
& \multicolumn{3}{c}{{\tt\scriptsize rs228}} && \multicolumn{3}{c}{{\tt\scriptsize rs365}} \\
\cline{2-4} \cline{6-8}\\[-0.5em]
& Time (s) & \# Iter. & Objective &
& Time (s) & \# Iter. & Objective  \\
\cline{2-4} \cline{6-8}\\[-0.5em]
SeDuMi (high) & 1\,655 & 21  & 64.71  &
	   &  *** & ***  & ***  \\
SeDuMi (low) & 809 & 10 & 64.80 &	
					& *** & *** & *** \\
SCS (direct) & 2\,338 & 2\,000 &  62.06 &
	& 34\,497  & 2\,000  & 44.02    \\
{{\CDCS}}-primal & 94 & 400  & 64.65 & &  321  & 401  & 63.37  \\
{{\CDCS}}-dual & 84  & 341 & 64.76 & & 240  & 265   & 63.69  \\
{{\CDCS}}-hsde & 79  & 361 & 64.87 & & 332 & 442   & 63.64\\
\midrule
& \multicolumn{3}{c}{{\tt\scriptsize rs1555}} && \multicolumn{3}{c}{{\tt\scriptsize rs1907}} \\
\cline{2-4} \cline{6-8}\\[-0.5em]
& Time (s) & \# Iter. & Objective &
& Time (s) & \# Iter. & Objective  \\
\cline{2-4} \cline{6-8}\\[-0.5em]
SeDuMi (high)& *** & *** & ***   &
	  & *** & ***  & *** \\
SeDuMi (low) & ***  & *** & ***   &
					& ***  & *** & ***  \\
SCS (direct) & 139\,314   & 2\,000 & 34.20 &
	& 50\,047  & 2\,000 & 45.89  \\
{{\CDCS}}-primal & 1\,721 & 2\,000 & 61.22 &
			 &  330  & 349 & 62.87  \\
{{\CDCS}}-dual & 317  & 317 & 69.54  &	& 271  & 252  & 63.30    \\
{{\CDCS}}-hsde & 1\,413  & 2\,000 & 61.36  &	& 393 & 414 & 63.14 \\
\bottomrule[1pt]
\end{tabular*}
\scriptsize
\raggedright
\end{table*}

Total CPU time, number of iterations, and terminal objective values are presented in {Table~\ref{T:ResultsNonChordal}}. For all problems, the algorithms in {{\CDCS}} (primal, dual and hsde) are all much faster than either SCS or SeDuMi. In addition, SCS never terminates succesfully, while the objective value returned by {{\CDCS}} is always within 2\% of the high-accuracy solutions returned by SeDuMi (when this could be computed).  {The residuals listed in Tables~\ref{T:ResultsErrors} and~\ref{T:SedumiSCSresiduals} suggest that this performance difference might be due to poor problem scaling in SCS.}

{The advantages of the algorithms proposed in this work are evident from Table~\ref{T:ResultsNonChordalIteration}: the average CPU time per iteration in {{\CDCS}}-hsde is approximately  $22 \times$, $24 \times$, $ 28 \times$, and $ 105 \times$ faster compared to SCS for problems \probname{rs200}, \probname{rs365}, \probname{rs1907}, and \probname{rs1555}, respectively.}
{The results for average CPU time per iteration also demonstrate that the computational complexity of all three algorithms in {{\CDCS}} (primal, dual, and hsde) is independent of the original problem size:} problems  \probname{rs35} and  \probname{rs228} have similar cone size $n$ and the same number of constraints $m$, yet the average CPU time for the latter is approximately 5$\times$ smaller. {This can be explained by noticing that for all test problems considered here the number of constraints $m$ is moderate,} so the overall complexity of our algorithms is dominated by the conic projection. As stated in Proposition~\ref{T:FlopsConic}, this depends only on the size and number of the maximal cliques, not on the size of the original PSD cone. A more detailed investigation of how the number of maximal cliques, their size, and the number of constraints affect the performance of {\CDCS} is presented next.

\begin{table}
	\centering
	\renewcommand\arraystretch{0.9}
	\caption{{Average CPU time per iteration (in seconds) for the nonchordal SDPs form~\cite{andersen2010implementation}.}}
	\label{T:ResultsNonChordalIteration}
	\begin{tabular*}{\linewidth}{@{\extracolsep{\fill}}r rrrrrr}
		\toprule[1pt]
		& {\tt\scriptsize rs35} & {\tt\scriptsize rs200}
		& {\tt\scriptsize rs228} & {\tt\scriptsize rs365}
		& {\tt\scriptsize rs1555} & {\tt\scriptsize rs1907}\\
		\cline{2-7}\\[-0.5em]
		SCS (direct)  & 1.188 & 4.847 & 1.169 & 17.250 & 69.590 & 25.240\\
		{{\CDCS}}-primal  & 0.944 & 0.258  & 0.224 & 0.715 &  0.828 & 0.833  \\
		{{\CDCS}}-dual &  1.064 & 0.263  &  0.232 & 0.774 & 0.791 & 0.920\\
		{{\CDCS}}-hsde &  1.005& 0.222 & 0.212 & 0.735  & 0.675  & 0.893\\
		\bottomrule[1pt]
	\end{tabular*}
\end{table}

\subsection{Random SDPs with block-arrow patterns} \label{Section:BlockArrow}

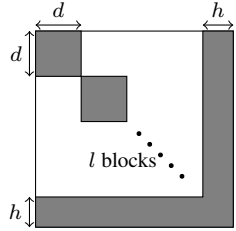
\begin{figure}
\centering
\setlength{\abovecaptionskip}{0pt}
\setlength{\belowcaptionskip}{1em}
\tikzset{decorate sep/.style 2 args=
{decorate,decoration={shape backgrounds,shape=circle,shape size=#1,shape sep=#2}}}
\begin{tikzpicture}[scale=0.8]
	 \draw (0,0) rectangle  (3.25,3.25);
	 \draw[fill=black!50] (0,3.25) rectangle  (0.75,2.5);
	 \draw[fill=black!50] (0.75,2.5) rectangle  (1.5,1.75);
	 \draw[fill=black!50] (0,0)--(3.25,0)--(3.25,3.25)--
	 				      (2.75,3.25)--(2.75,0.5)--(0,0.5)--(0,0);
	 \draw[decorate sep={0.5mm}{2mm},fill] (1.7,1.55)--
	 					   node[left] {\small $l$ blocks }(2.55,0.7);
	 \draw[<->] (0,3.35) -- node[above] {\small $d$} (0.75,3.35) ;
	 \draw[<->] (-0.1,2.5) -- node[left] {\small $d$} (-0.1,3.25) ;
	 \draw[<->] (2.75,3.35) -- node[above] {\small $h$} (3.25,3.35) ;
	 \draw[<->] (-0.1,0) -- node[left] {\small $h$} (-0.1,0.5) ;
\end{tikzpicture}
\caption{Block-arrow sparsity pattern (dots indicate repeating diagonal blocks). The parameters are: the number of blocks, $l$; block size, $d$; the width of the arrow head, $h$.}
\label{F:BlockArrowSDP}
\end{figure}

{To examine the influence of the number of maximal cliques, their size, and the number of constraints on the computational cost of Algorithms~\ref{A:ADMMPrimal}--\ref{A:ADMMhsde}}, we considered randomly generated SDPs with a ``block-arrow'' aggregate sparsity pattern, {illustrated} in {Fig.}~\ref{F:BlockArrowSDP}. Such a sparsity pattern is characterized by: the number of blocks, $l$; the block size, $d$; and the size of the arrow head, $h$. {The associated PSD cone has dimension $ld+h$.} The block-arrow sparsity pattern is chordal, with $l$ maximal cliques all of the same size $d+h$. The effect of the number of constraints in the SDP, $m$, is investigated as well, and numerical results are presented below for the following scenarios:
\begin{enumerate}
  \item {Fix $l=100$, $d=10$, $h=20$, and vary the number of constraints, $m$;}
  \item {Fix $m=200$, $d=10$, $h=20$, and vary $l$ (hence, the number of maximal cliques);}
  \item {Fix $m = 200$, $l=50$, $h=10$, and vary $d$ (hence, the size of {the}  maximal cliques).}
\end{enumerate}

In our computations, the problem data are generated randomly using the following procedure. {First, we generate random symmetric matrices $A_1,\,\ldots,\,A_m$ with block-arrow sparsity pattern, whose nonzero entries are drawn from the uniform distribution $U(0,1)$ on the open interval $(0,1)$. Second,} a strictly primal feasible matrix $X_{\text{f}}\in \mathbb{S}_+^n(\mathcal{E},0)$ is constructed as $X_{\text{f}} = W + \alpha I$, where $W \in \mathbb{S}^n(\mathcal{E},0)$ is randomly generated with entries from $U(0,1)$ and $\alpha$ is chosen to guarantee $X_{\text{f}} \succ 0$. The vector $b$ in the primal {equality constraints} is {then} computed {such that} $b_i = \langle A_i,X_{\text{f}} \rangle$ for all $i = 1, \ldots, m$. {Finally,} the matrix $C$ in the dual constraint is constructed as $C = Z_{\text{f}} + \sum_{i=1}^my_iA_i$, where $y_1,\,\ldots,\,y_m$ are drawn from $U(0,1)$ and $Z_{\text{f}} \succ 0$ is generated similarly to $X_{\text{f}}$.

\begin{figure*}
    \centering
	\includegraphics[width=0.9\textwidth]{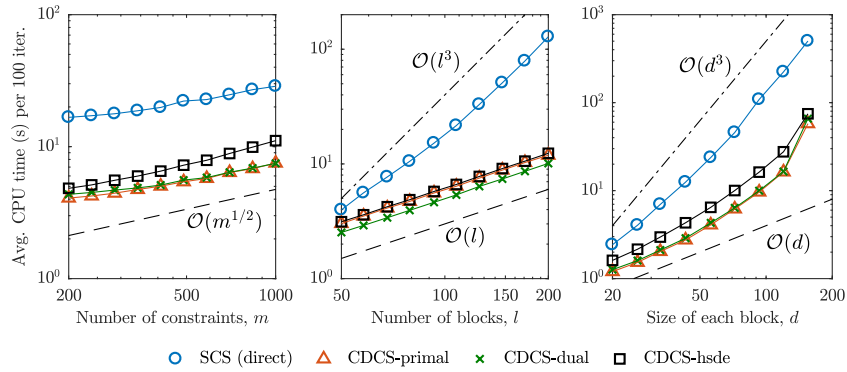}
    \caption{Average CPU time (in seconds) per 100 iterations for SDPs with block-arrow patterns. Left to right: varying the number of constraints; varying the number of blocks; varying the block size.   }
    \label{F:BlcokArrowResult}
\end{figure*}

\begin{table}
    \centering
    \caption{Average CPU time ($\times 10^{-2}$ s) required by the affine projection steps in {\CDCS}-primal, {\CDCS}-dual, and {\CDCS}-hsde as a function of the number of constraints ($m$) for $l=100$, $d=10$, and $h=20$.}
    \label{T:ResultsLinProj}
    \begin{tabular*}{\linewidth}{@{\extracolsep{\fill}}r| rrrrrrrrrr}
    \toprule[1pt]
     $m$ & 200 & 239 & 286 & 342  & 409 & 489 & 585 & 699 & 836 & 1000 \\
    \midrule 
    {\CDCS}-primal  & 1.05 & 1.21  & 1.40 & 1.63 &  1.90 & 2.22 & 2.60 & 3.12 & 3.59 & 4.29 \\
    {\CDCS}-dual &  1.10 & 1.26  &  1.46 & 1.67 & 1.94 & 2.28 & 2.65  & 3.16 &  3.66 & 4.31\\
    {\CDCS}-hsde &  1.84 & 2.14 & 2.55 & 2.95  & 3.50  & 4.12 & 4.85  & 5.80 & 6.81 &  8.04\\
    \bottomrule[1pt]
    \end{tabular*}
\end{table}

The average CPU time per 100 iterations for the first-order solvers is {plotted} in Figure~\ref{F:BlcokArrowResult}. As already observed in the previous sections, in all three test scenarios the algorithms in {{\CDCS}} are faster than SCS, when the latter is {used to solve the original SDPs (before chordal decomposition)}. Of course, as one would expect, {the} computational cost grows when either the number of constraints, the size of the maximal cliques, or their number is increased. {Note, however, that the CPU time per iteration of CDCS grows more slowly than that of SCS as a function of the number of maximal cliques, which is the benefit of considering smaller PSD cones in {{\CDCS}}.} Precisely, the CPU time per iteration of CDCS increases linearly when the number of cliques $l$ is raised, as expected from Proposition~\ref{T:FlopsConic}; instead, the CPU time per iteration of SCS grows cubically, since the eigenvalue decomposition on the original cone requires $\mathcal{O}(l^3)$ flops (note that when $d$ and $h$ are fixed, $(ld+h)^3 = \mathcal{O}(l^3)$). Finally, the results in Table~\ref{T:ResultsLinProj} confirm the analysis in Propositions~\ref{T:FlopsLinP&D} and~\ref{T:FlopsLinHsde_final}, according to which the CPU time required in the affine projection of {{\CDCS}}-hsde was approximately twice larger than that of {{\CDCS}}-primal or {{\CDCS}}-dual. {On the other hand, the increase in computational cost with the number of constraints $m$ is slower than predicted by Propositions~\ref{T:FlopsLinP&D} and~\ref{T:FlopsLinHsde_final} due to the fact that, contrary to the complexity analysis presented in Section~\ref{Section:Complexity}, our implementation of Algorithms~\ref{A:ADMMPrimal}--\ref{A:ADMMhsde} takes advantage of sparse matrix operations where possible.}

\section{Conclusion} \label{sec:conclusion}

In this paper, we have presented a conversion framework for large-scale SDPs characterized by chordal sparsity. This framework is analogous to the conversion techniques for IPMs of~\cite{fukuda2001exploiting, kim2011exploiting}, but is more suitable for the application of FOMs. We have then developed efficient ADMM algorithms for sparse SDPs in either primal or dual standard form, and for their homogeneous self-dual embedding. In all cases, a single iteration of our ADMM algorithms only requires parallel projections onto small PSD cones and a projection onto an affine subspace, both of which can be carried out efficiently. In particular, when the number of constraints $m$ is moderate the complexity of each iteration is determined by the size of the largest maximal clique, not the size of the original problem. This enables us to solve large, sparse conic problems that are beyond the reach of standard interior-point and/or other first-order methods.

All our algorithms have been made available in the open-source MATLAB solver {{\CDCS}}. Numerical simulations on benchmark problems, including selected sparse problems from SDPLIB, large and sparse SDPs with a nonchordal sparsity pattern, and SDPs with a block-arrow sparsity pattern, demonstrate that our methods can significantly reduce the total CPU time requirement compared to the state-of-the-art interior-point solver SeDuMi \cite{sturm1999using} and the efficient first-order solver SCS~\cite{scs}. We remark that the current implementation of our algorithms is sequential, but many steps can be carried out in parallel, so further computational gains may be achieved by taking full advantage of distributed computing architectures. Besides, it would be interesting to integrate some acceleration techniques (\emph{e.g.},~\cite{themelis2016supermann,falt2016line}) that {promise} to improve the convergence performance of ADMM in practice.

Finally, we note that the conversion framework we have proposed relies on chordal sparsity, but there exist large SDPs which do not have this property. An example with applications in many areas {is that of} SDPs from sum-of-squares relaxations of polynomial optimization problems. Future work should therefore explore whether and {to which extent} first order methods can be used to take advantage other types of sparsity and structure.

{
\begin{acknowledgements}
    The authors would like to thank the Associate Editor and the anonymous reviewers, whose invaluable comments contributed to improving the quality of our original manuscript.
\end{acknowledgements}

\appendix
\section*{Appendix}

\textbf{A\hspace{5pt}Proof of Proposition~\ref{T:FlopsLinP&D}.}
Since~\eqref{E:OptCondMinYPrimal} and~\eqref{E:OptCondMinYDual} are the same modulo scaling, we only consider the former. Also, we drop the superscript $(n)$ to lighten the notation. Recall that $H_k^Tx_k$ is an indexing operation and requires no flops, and let
    \begin{equation}
    \label{e:bnflopsdef}
        \hat{b} := \sum_{k=1}^{p} H_k^T\left( x_k+\rho^{-1}\lambda_k\right) - \rho^{-1}c  \in \mathbb{R}^{n^2}.
   \end{equation}
    After a suitable block elimination and writing $AD^{-1}A^T = LL^T$, the solution of~\eqref{E:OptCondMinYPrimal} is given by
    \begin{subequations}  \label{E:KKTflops}
        \begin{align}
	        LL^T y &= AD^{-1}\hat{b} - b,  \label{E:KKTflops_s2}\\
            x &= D^{-1} \left(\hat{b} - A^Ty \right).
            \label{E:KKTflops_s1}
        \end{align}
    \end{subequations}
    Computing $x$ and $y$ cost $(4m+p+3)n^2+2m^2+2n_d$ flops, counted as the sum of:

    \begin{enumerate}[(i)]
          \item  $(p+1)n^2 + 2n_d$ flops to form $\hat{b}$:  no flops to multiply by $H_k$, $2|\mathcal{C}_k|^2$ flops to compute $x_k+\rho^{-1}\lambda_k$, $n^2$ flops to calculate $\rho^{-1}c$, and $(p-1)n^2+n^2$ flops to sum all addends in~\eqref{e:bnflopsdef}.
          \item $(2m+1)n^2$ flops to compute $AD^{-1}\hat{b} - b$: $n^2$ flops to compute $D^{-1}\hat{b}$ since $D$ is diagonal, $(2n^2-1)m$ flops to multiply by $A$, and $m$ flops to subtract $b$.
          \item $2m^2$ flops to compute $y$ via forward and backward substitutions using~\eqref{E:KKTflops_s2}.
          \item  $(2m + 1)n^2$ flops to compute $x$ via~\eqref{E:KKTflops_s1}: $(2m-1)n^2$ flops to find $A^T y$, $n^2$ flops to subtract it from $\hat{b}$, and $n^2$ flops to multiply by $D^{-1}$.
    \end{enumerate}

\vspace{1em}
\noindent
\textbf{B\hspace{5pt}Proof of Proposition~\ref{T:FlopsLinHsde_final}.}  Consider the ``inner'' system~\eqref{E:SecondBlockElimination} first. Partition the vectors $\sigma_1$ and $\sigma_2$ as
    $$
        \sigma_1 = \begin{bmatrix} \sigma_{11} \\ \sigma_{12} \end{bmatrix}, \qquad
        \sigma_2 = \begin{bmatrix} \sigma_{21} \\ \sigma_{22} \end{bmatrix},
    $$
    where $\sigma_{11} \in \mathbb{R}^{n^2}$, $\sigma_{12},\,\sigma_{22} \in \mathbb{R}^{n_d}$, and $\sigma_{21} \in \mathbb{R}^{m}$.  The vectors $\nu_1$ and $\nu_2$ on the right-hand side of~\eqref{E:SecondBlockElimination} can be partitioned in a similar way.
    Recalling the definition of the matrix $\hat{A}$ from~\eqref{E:MDef},~\eqref{E:SecBlkEliResult} becomes
    \begin{equation}
    \label{e:AppendixEqSigma2}
        \begin{bmatrix} \sigma_{21} \\ \sigma_{22} \end{bmatrix} =
        \begin{bmatrix} \nu_{21} - A\sigma_{11} \\ \nu_{22} - H\sigma_{11}+\sigma_{12}  \end{bmatrix}.
    \end{equation}
    To calculate $\sigma_{11}$ and $\sigma_{12}$ one needs to solve~\eqref{E:SecBlkEliResult_s1}, which after partitioning all variables can be rewritten as
    \begin{equation} \label{E:SecBlkEliFlops}
        \begin{bmatrix} \left(I+D + A^TA \right) & -H^T \\ -H & 2I\end{bmatrix} \begin{bmatrix} \sigma_{11} \\ \sigma_{12} \end{bmatrix} =
		\begin{bmatrix}
		\nu_{11} + A^T\nu_{21}+H^T\nu_{22} \\
		\nu_{12} - \nu_{22}
		\end{bmatrix}.
    \end{equation}
    Eliminating $\sigma_{12}$ from the first block equation results in
	\begin{subequations}
		\begin{align}
		\left(I + \frac{1}{2}D + A^TA \right)\sigma_{11} &= \nu_{11} + A^T\nu_{21}+ \frac{1}{2}H^T\left( \nu_{12} + \nu_{22} \right), \label{E:ThirBlkFlops_s1}\\
		\sigma_{12}&=\frac{1}{2}(\nu_{12} - \nu_{22} + H\sigma_{11}) \label{E:ThirBlkFlops_s2}.
		\end{align}
	\end{subequations}
    After defining $P:=I+\frac{1}{2}D$ and
    $\eta := \nu_{11} + A^T\nu_{21}+ \frac{1}{2}H^T\left( \nu_{12} + \nu_{22} \right)$
    to lighten the notation, an application of the matrix inversion lemma to~\eqref{E:ThirBlkFlops_s1} yields
    \begin{equation} \label{E:MatrixInvFlops}
        \sigma_{11} = P^{-1}\eta - P^{-1}A^T(I+AP^{-1}A^T)^{-1}AP^{-1}\eta.
    \end{equation}
    We are now in a position to count the flops required to solve the ``inner'' linear system. First, computing $\sigma_{11}$ via~\eqref{E:MatrixInvFlops} requires a total $(6m+p+3)n^2+2m^2-m$ flops, counted as follows:

	\begin{enumerate}[(i)]
	  \item $(2m+p+1)n^2$ flops to form $\eta$;
      \item $n^2$ flops to compute $P^{-1}\eta$, since $P$ is an $n^2 \times n^2$ diagonal matrix;
      \item $(2n^2-1)m$ flops to calculate $AP^{-1}\eta$;
      \item $2m^2$ flops to form the vector $(I+AP^{-1}A^T)^{-1}AP^{-1}\eta$ using forward and backward substitutions (we assume that the Cholesky decomposition $I+AP^{-1}A^T = LL^T$ has been cached);
      \item $(2m-1)n^2$ flops to find $A^T(I+AP^{-1}A^T)^{-1}AP^{-1}\eta$;
      \item {$2n^2$} flops to compute $ \sigma_{11}$ via~\eqref{E:MatrixInvFlops} given $P^{-1}\eta$ and $A^T(I+AP^{-1}A^T)^{-1}AP^{-1}\eta$.
    \end{enumerate}

	\noindent
	Once $\sigma_{11}$ is known, $\sigma_{12}$ is found from~\eqref{E:ThirBlkFlops_s2} with $3 n_d$ flops because the product $H \sigma_{11}$ is simply an indexing operation and costs no flops. Given $\sigma_{11}$ and $\sigma_{12}$, computing $\sigma_{21}$ and $\sigma_{22}$ from~\eqref{e:AppendixEqSigma2} requires $2m n^2 + 2n_d$ flops, so the ``inner'' linear system~\eqref{E:SecondBlockElimination} costs a total of $(8m+2p+3)n^2+2m^2 - m +5n_d$ flops.

After the inner system has been solved, we see that computing $\hat{u}_1$ from~\eqref{E:MatrixInverseResult} requires $(8m+2p+9)n^2+2m^2 + 5m +17n_d-1$ flops in total:

    \begin{enumerate}[(i)]
      \item $2(n^2+2n_d+m)$ flops to compute $\omega_1 - \omega_2\zeta$;
      \item $(8m+2p+3)n^2+2m^2 - m +5n_d$ flops to solve the ``inner'' linear system $M^{-1}(\omega_1 - \omega_2\zeta)$;
      \item $2(n^2+2n_d+m)-1$ flops to compute $\zeta^TM^{-1}(\omega_1 - \omega_2\zeta) \in \mathbb{R}$;
      \item $n^2+2n_d+m$ flops to calculate $\hat{\zeta} \cdot \zeta^TM^{-1}(\omega_1 - \omega_2\zeta)$;
      \item $n^2+2n_d+m$ flops to compute $\hat{u}_1 = M^{-1}(\omega_1 - \omega_2\zeta) - \hat{\zeta} \cdot \zeta^TM^{-1}(\omega_1 - \omega_2\zeta)$.
    \end{enumerate}

	\noindent
	Summing this to the $2(n^2+2n_d+m)$ flops required to calculate $\hat{u}_2$ using~\eqref{E:FirstBlockElimination} yields the desired result.

\vspace{1em}
\noindent
\textbf{C\hspace{5pt}Proof of Proposition~\ref{T:FlopsConic}.}
The conic projection~\eqref{E:XkUpdate} in Algorithm~\ref{A:ADMMPrimal}  amounts to projecting the matrices
$$\mat\left( H_k x^{(n+1)} -  {\rho}^{-1}\lambda_k^{(n)}\right) \in \mathbb{S}^{|\mathcal{C}_k|}, \quad k = 1, \ldots, p$$
onto the PSD cone $\mathbb{S}^{|\mathcal{C}_k|}_+$. Computing $H_k x^{(n+1)} -  {\rho}^{-1}\lambda_k^{(n)}$ requires $2|\mathcal{C}_k|^2$ flops, while a PSD projection using a full eigenvalue decomposition costs $\mathcal{O}(|\mathcal{C}_k|^3)$ flops to leading order, so the overall number of flops is $\mathcal{O}(\sum_{k=1}^p|\mathcal{C}_k|^3)$. The same argument holds for the conic projection~\eqref{E:zkUpdate} in Algorithm~\ref{A:ADMMDual}.

In Algorithm~\ref{A:ADMMhsde}, instead, the projection is onto the cone
$
\mathcal{K} := \mathbb{R}^{n^2} \times \mathcal{S} \times \mathbb{R}^{m} \times \mathbb{R}^{n_d} \times \mathbb{R}_{+}.
$
Nothing needs to be done to project onto $\mathbb{R}^{n^2}$, $\mathbb{R}^{m}$ and $\mathbb{R}^{n_d}$, while the projection of $a\in\mathbb{R}$ onto $\mathbb{R}_{+}$ is given by $\max\{0,a\}$ and requires no flops according to our definition. Finally, projecting onto $\mathcal{S}$ requires eigenvalue decompositions of the matrices $\mat(x_k)$, $k=1,\ldots,p$, with a leading-order cost of $\mathcal{O}(\sum_{k=1}^p|\mathcal{C}_k|^3)$ flops.
}

\bibliographystyle{spmpsci}      
\bibliography{Reference}   

\begin{thebibliography}{10}
\providecommand{\url}[1]{{#1}}
\providecommand{\urlprefix}{URL }
\expandafter\ifx\csname urlstyle\endcsname\relax
  \providecommand{\doi}[1]{DOI~\discretionary{}{}{}#1}\else
  \providecommand{\doi}{DOI~\discretionary{}{}{}\begingroup
  \urlstyle{rm}\Url}\fi

\bibitem{agler1988positive}
Agler, J., Helton, W., McCullough, S., Rodman, L.: Positive semidefinite
  matrices with a given sparsity pattern.
\newblock Linear Algebra Appl. \textbf{107}, 101--149 (1988)

\bibitem{alizadeh1998primal}
Alizadeh, F., Haeberly, J.P.A., Overton, M.L.: Primal-dual interior-point
  methods for semidefinite programming: convergence rates, stability and
  numerical results.
\newblock SIAM J. Optim. \textbf{8}(3), 746--768 (1998)

\bibitem{andersen2011interior}
Andersen, M., Dahl, J., Liu, Z., Vandenberghe, L.: Interior-point methods for
  large-scale cone programming.
\newblock In: Optimization for machine learning, pp. 55--83. MIT Press (2011)

\bibitem{andersen2010implementation}
Andersen, M.S., Dahl, J., Vandenberghe, L.: Implementation of nonsymmetric
  interior-point methods for linear optimization over sparse matrix cones.
\newblock Math. Program. Comput. \textbf{2}(3-4), 167--201 (2010)

\bibitem{BGSBfeasibility2017}
Banjac, G., Goulart, P., Stellato, B., Boyd, S.: Infeasibility detection in the
  alternating direction method of multipliers for convex optimization.
\newblock optimization-online.org  (2017).
\newblock
  \urlprefix\url{http://www.optimization-online.org/DB_HTML/2017/06/6058.html}

\bibitem{blair1993introduction}
Blair, J.R., Peyton, B.: An introduction to chordal graphs and clique trees.
\newblock In: Graph theory and sparse matrix computation, pp. 1--29. Springer
  (1993)

\bibitem{borchers1999sdplib}
Borchers, B.: {SDPLIB} 1.2, a library of semidefinite programming test
  problems.
\newblock Optim. Methods Softw. \textbf{11}(1-4), 683--690 (1999)

\bibitem{boyd1994linear}
Boyd, S., El~Ghaoui, L., Feron, E., Balakrishnan, V.: Linear Matrix
  Inequalities in System and Control Theory.
\newblock SIAM (1994)

\bibitem{boyd2011distributed}
Boyd, S., Parikh, N., Chu, E., Peleato, B., Eckstein, J.: Distributed
  optimization and statistical learning via the alternating direction method of
  multipliers.
\newblock Found. Trends{\textregistered} Mach. Learn. \textbf{3}(1), 1--122
  (2011)

\bibitem{boyd2004convex}
Boyd, S., Vandenberghe, L.: Convex optimization.
\newblock Cambridge University Press (2004)

\bibitem{burer2003semidefinite}
Burer, S.: Semidefinite programming in the space of partial positive
  semidefinite matrices.
\newblock SIAM J. Optim. \textbf{14}(1), 139--172 (2003)

\bibitem{dall2013distributed}
Dall'Anese, E., Zhu, H., Giannakis, G.B.: Distributed optimal power flow for
  smart microgrids.
\newblock IEEE Trans. Smart Grid \textbf{4}(3), 1464--1475 (2013)

\bibitem{CSparse}
Davis, T.: Direct Methods for Sparse Linear Systems.
\newblock SIAM (2006)

\bibitem{davis2011university}
Davis, T.A., Hu, Y.: The {University of Florida} sparse matrix collection.
\newblock ACM Trans. Math. Softw. \textbf{38}(1), 1 (2011)

\bibitem{falt2016line}
F{\"a}lt, M., Giselsson, P.: Line search for generalized alternating
  projections.
\newblock arXiv preprint arXiv:1609.05920  (2016)

\bibitem{fujisawa2009user}
Fujisawa, K., Kim, S., Kojima, M., Okamoto, Y., Yamashita, M.: User's manual
  for {SparseCoLO}: Conversion methods for sparse conic-form linear
  optimization problems.
\newblock Tech. rep., Research Report B-453, Tokyo Institute of Technology,
  Tokyo 152-8552, Japan (2009)

\bibitem{fukuda2001exploiting}
Fukuda, M., Kojima, M., Murota, K., Nakata, K.: {Exploiting sparsity in
  semidefinite programming via matrix completion I: General framework}.
\newblock SIAM J. Optim. \textbf{11}(3), 647--674 (2001)

\bibitem{gabay1976dual}
Gabay, D., Mercier, B.: A dual algorithm for the solution of nonlinear
  variational problems via finite element approximation.
\newblock Comput. Math. Appl. \textbf{2}(1), 17--40 (1976)

\bibitem{ghadimi2015optimal}
Ghadimi, E., Teixeira, A., Shames, I., Johansson, M.: Optimal parameter
  selection for the alternating direction method of multipliers {(ADMM)}:
  quadratic problems.
\newblock IEEE Trans. Automat. Contr. \textbf{60}(3), 644--658 (2015)

\bibitem{glowinski1975approximation}
Glowinski, R., Marroco, A.: Sur l'approximation, par {\'e}l{\'e}ments finis
  d'ordre un, et la r{\'e}solution, par p{\'e}nalisation-dualit{\'e} d'une
  classe de probl{\`e}mes de dirichlet non lin{\'e}aires.
\newblock Revue fran{\c{c}}aise d'automatique, informatique, recherche
  op{\'e}rationnelle. Analyse num{\'e}rique \textbf{9}(2), 41--76 (1975)

\bibitem{godsil2013algebraic}
Godsil, C., Royle, G.F.: Algebraic graph theory.
\newblock Springer Science \& Business Media (2013)

\bibitem{griewank1984existence}
Griewank, A., Toint, P.L.: On the existence of convex decompositions of
  partially separable functions.
\newblock Math. Program. \textbf{28}(1), 25--49 (1984)

\bibitem{grone1984positive}
Grone, R., Johnson, C.R., S{\'a}, E.M., Wolkowicz, H.: Positive definite
  completions of partial hermitian matrices.
\newblock Linear Algebra Appl. \textbf{58}, 109--124 (1984)

\bibitem{helmberg1996interior}
Helmberg, C., Rendl, F., Vanderbei, R.J., Wolkowicz, H.: An interior-point
  method for semidefinite programming.
\newblock SIAM J. Optim. \textbf{6}(2), 342--361 (1996)

\bibitem{kakimura2010direct}
Kakimura, N.: A direct proof for the matrix decomposition of chordal-structured
  positive semidefinite matrices.
\newblock Linear Algebra Appl. \textbf{433}(4), 819--823 (2010)

\bibitem{Kalbat2015Fast}
Kalbat, A., Lavaei, J.: A fast distributed algorithm for decomposable
  semidefinite programs.
\newblock In: Proc. 54th IEEE Conf. Decis. Control, pp. 1742--1749 (2015)

\bibitem{kim2011exploiting}
Kim, S., Kojima, M., Mevissen, M., Yamashita, M.: Exploiting sparsity in linear
  and nonlinear matrix inequalities via positive semidefinite matrix
  completion.
\newblock Math. Program. \textbf{129}(1), 33--68 (2011)

\bibitem{liu2017new}
Liu, Y., Ryu, E.K., Yin, W.: A new use of douglas-rachford splitting and {ADMM}
  for identifying infeasible, unbounded, and pathological conic programs.
\newblock arXiv preprint arXiv:1706.02374  (2017)

\bibitem{lofberg2005yalmip}
Lofberg, J.: {YALMIP}: A toolbox for modeling and optimization in {MATLAB}.
\newblock In: IEEE Int. Symp. Comput. Aided Control Sys. Des., pp. 284--289.
  IEEE (2004)

\bibitem{Madani2015ADMM}
Madani, R., Kalbat, A., Lavaei, J.: {ADMM} for sparse semidefinite programming
  with applications to optimal power flow problem.
\newblock In: Proc. 54th IEEE Conf. Decis. Control, pp. 5932--5939 (2015)

\bibitem{malick2009regularization}
Malick, J., Povh, J., Rendl, F., Wiegele, A.: Regularization methods for
  semidefinite programming.
\newblock SIAM J. Optim. \textbf{20}(1), 336--356 (2009)

\bibitem{ODonoghue2016}
O'Donoghue, B., Chu, E., Parikh, N., Boyd, S.: Conic optimization via operator
  splitting and homogeneous self-dual embedding.
\newblock J. Optim. Theory Appl. \textbf{169}(3), 1042--1068 (2016)

\bibitem{scs}
O'Donoghue, B., Chu, E., Parikh, N., Boyd, S.: {SCS}: Splitting conic solver,
  version 1.2.6.
\newblock \url{https://github.com/cvxgrp/scs} (2016)

\bibitem{papachristodoulou2013sostools}
Papachristodoulou, A., Anderson, J., Valmorbida, G., Prajna, S., Seiler, P.,
  Parrilo, P.: {SOSTOOLS} version 3.00 sum of squares optimization toolbox for
  {MATLAB}.
\newblock arXiv preprint arXiv:1310.4716  (2013)

\bibitem{raghunathan2014alternating}
Raghunathan, A.U., Di~Cairano, S.: Alternating direction method of multipliers
  for strictly convex quadratic programs: Optimal parameter selection.
\newblock In: Proc. American Control Conf., pp. 4324--4329. IEEE (2014)

\bibitem{saad2003iterative}
Saad, Y.: Iterative methods for sparse linear systems.
\newblock SIAM (2003)

\bibitem{sturm1999using}
Sturm, J.F.: Using {SeDuMi} 1.02, a {MATLAB} toolbox for optimization over
  symmetric cones.
\newblock Optim. Methods Softw. \textbf{11}(1-4), 625--653 (1999)

\bibitem{sun2014decomposition}
Sun, Y., Andersen, M.S., Vandenberghe, L.: Decomposition in conic optimization
  with partially separable structure.
\newblock SIAM J. Optim. \textbf{24}(2), 873--897 (2014)

\bibitem{sun2015decomposition}
Sun, Y., Vandenberghe, L.: Decomposition methods for sparse matrix nearness
  problems.
\newblock SIAM J. Matrix Anal. Appl. \textbf{36}(4), 1691--1717 (2015)

\bibitem{tarjan1984simple}
Tarjan, R.E., Yannakakis, M.: Simple linear-time algorithms to test chordality
  of graphs, test acyclicity of hypergraphs, and selectively reduce acyclic
  hypergraphs.
\newblock SIAM J. Comput. \textbf{13}(3), 566--579 (1984)

\bibitem{themelis2016supermann}
Themelis, A., Patrinos, P.: {SuperMann}: a superlinearly convergent algorithm
  for finding fixed points of nonexpansive operators.
\newblock arXiv preprint arXiv:1609.06955  (2016)

\bibitem{vandenberghe2014chordal}
Vandenberghe, L., Andersen, M.S.: Chordal graphs and semidefinite optimization.
\newblock Found. Trends{\textregistered} Optim. \textbf{1}(4), 241--433 (2014)

\bibitem{vandenberghe1996semidefinite}
Vandenberghe, L., Boyd, S.: Semidefinite programming.
\newblock SIAM Review \textbf{38}(1), 49--95 (1996)

\bibitem{wen2010alternating}
Wen, Z., Goldfarb, D., Yin, W.: Alternating direction augmented lagrangian
  methods for semidefinite programming.
\newblock Math. Program. Comput. \textbf{2}(3-4), 203--230 (2010)

\bibitem{yan2016self}
Yan, M., Yin, W.: Self equivalence of the alternating direction method of
  multipliers.
\newblock In: Splitting Methods in Communication, Imaging, Science, and
  Engineering, pp. 165--194. Springer (2016)

\bibitem{yannakakis1981computing}
Yannakakis, M.: Computing the minimum fill-in is \protect{NP-complete}.
\newblock SIAM J. Algebraic Discrete Methods \textbf{2}, 77--79 (1981)

\bibitem{ye2011interior}
Ye, Y.: Interior point algorithms: theory and analysis.
\newblock John Wiley \& Sons (2011)

\bibitem{ye1994nl}
Ye, Y., Todd, M.J., Mizuno, S.: An $\mathcal{O}\sqrt n l $-iteration
  homogeneous and self-dual linear programming algorithm.
\newblock Math. Oper. Res. \textbf{19}(1), 53--67 (1994)

\bibitem{zhao2010newton}
Zhao, X.Y., Sun, D., Toh, K.C.: A newton-cg augmented lagrangian method for
  semidefinite programming.
\newblock SIAM J. Optim. \textbf{20}(4), 1737--1765 (2010)

\bibitem{zheng2017Exploiting}
Zheng, Y., Fantuzzi, G., Papachristodoulou, A.: Exploiting sparsity in the
  coefficient matching conditions in {Sum-of-Squares} programming using {ADMM}.
\newblock IEEE Control Systems Letters \textbf{1}(1), 80--85 (2017)

\bibitem{ZFPGWhsde2016}
Zheng, Y., Fantuzzi, G., Papachristodoulou, A., Goulart, P., Wynn, A.: Fast
  {ADMM} for homogeneous self-dual embedding of sparse {SDPs}.
\newblock IFAC PapersOnLine \textbf{50}(1), 8411--8416 (2017)

\bibitem{ZFPGWpd2016}
Zheng, Y., Fantuzzi, G., Papachristodoulou, A., Goulart, P., Wynn, A.: Fast
  {ADMM} for semidefinite programs with chordal sparsity.
\newblock In: Proc. American Control Conf., pp. 3335--3340. IEEE (2017)

\end{thebibliography}

\end{document}